\documentclass[referee,
sn-standardnature]{sn-jnl}

\usepackage{subcaption,array,tabularx,multirow,multicol} 
\jyear{2024}%

\usepackage{soul}
\usepackage{booktabs}
\begin{document}

\title[Escape from a well]{Escape from a potential well under forcing and damping with application to ship capsize}


\author*[1,2]{\fnm{Alex} \sur{McSweeney-Davis}}\email{alex@mcsweeney-davis.com}

\author[1]{\fnm{R.S.} \sur{MacKay}}\email{R.S.MacKay@warwick.ac.uk}

\author[1,3]{\fnm{Shibabrat} \sur{Naik}}\email{shibabratnaik@gmail.com}


\affil*[1]{\orgdiv{Mathematics Institute}, \orgname{University of Warwick}, 
\\ \city{Coventry}, \postcode{CV4 7AL}, 
\country{U.K.}}

\affil[2]{\footnotemark[2] \orgdiv{Equipe IDEFIX}, \orgname{INRIA Saclay} , \orgaddress{UMA, ENSTA Paris, 828, Boulevard des Marechaux}, \postcode{91762} \city{Palaiseau}, \country{France}}

\affil[3]{\footnote[2]{(Current address)} \orgdiv{Department of Mathematics}, \orgname{Hampton University}, \\ \orgaddress{\street{609 Norma B Harvey Road}, \city{Hampton}, \postcode{23669}, \state{Virginia}, \country{U.S.}}}

\abstract{Escape from a potential well occurs in a wide variety of physical systems from chemical reactions to ship capsize. In these situations escape often occurs by passage ``over'' a normally hyperbolic submanifold (NHS). This paper describes the computational implementation of an algorithm to identify the NHS of a two degree of freedom model for ship motion under the influence of damping and aperiodic forcing. Additionally we demonstrate how the stable manifolds of these submanifolds can be used to classify initial ship states as safe or unsafe.}




\maketitle
\tableofcontents
\clearpage
\markboth{\textit{Escape from a well}}{\textit{Escape from a well}}
\section{Introduction}
\label{intro}
Escape from a potential well is a phenomenon with widespread relevance, such as ionization of a hydrogen atom under an electromagnetic field in atomic physics~\cite{jaffe1999transition}, transport of defects in solid-state and semiconductor physics~\cite{eckhardt1995transition}, isomerization of clusters~\cite{komatsuzaki2001dynamical}, reaction rates in chemical physics~\cite{eyring1935activated,wigner1938transition,wiggins2001impenetrable}, buckling modes in structural mechanics~\cite{collins2012isomerization,zhong2018tube}, ship motion and capsize~\cite{virgin1989approx,thompson1996suppression,naik2017geometry}, escape and recapture of comets and asteroids in celestial mechanics~\cite{jaffe2002statistical,dellnitz2005transport}, and escape into inflation or recollapse to singularity in cosmology~\cite{deoliveira2002homoclinic}.

Although the autonomous Hamiltonian case and its extension to periodic forcing and damping are well-studied, developing the theory to allow aperiodic forcing is important for many real-world applications. There is substantial work on stochastic forcing, but the emphasis tends to be on probabilistic results with respect to the distribution of the forcing, whereas we want results for specified realisations of the forcing in order to devise safety criteria or control strategies. Even though the perturbation methods developed in~\cite{nayfeh1973nonlinear} are useful for deriving algebraic conditions for dynamical behaviours, a geometric viewpoint highlights the dynamical (phase space) mechanism of capsize. This was developed in studies by~\cite{thompson1996suppression,falzarano1992application} using the paradigm of the escape from a potential well for capsize and lobe dynamics as the mechanism of eroding the non-capsize region in the well. Approaches to include stochastic forcing are to be found in~\cite{naik2017geometry}, and dissipative, gyroscopic forcing in~\cite{zhong2020geometry}.  The dynamical systems mechanism of crossing a saddle in the presence of dissipation and forcing has been studied extensively in the context of chemical reactions to compute reaction rates accurately based on the dividing manifold~\cite{hernandez2010transition,pollak2005reaction,craven2017lagrangian}. 

The paper~\cite{bujorianu2021new} proposed a mathematical approach to quantifying escape from a potential well for systems with aperiodic time-dependent forcing and damping in the particular context of ship capsize. In this approach, the ship's state (configuration and velocities) is modelled by coupled non-linear non-autonomous dissipative ordinary differential equations. A criterion for capsize under aperiodic forcing was established based on the geometric viewpoint of dynamical systems theory. Here, we report on the computational implementation of the approach. 

As is standard for a non-autonomous system on a manifold $\mathbb{M}$, we consider the flow on the extended state space $\mathbb{M} \times \mathbb{R}$ (extended by time). In ~\cite{bujorianu2021new} the dynamics on the extended state space is studied with the help of a normally hyperbolic submanifold (NHS).
It is a (locally) invariant set that hovers on the brink of transition. It can be spanned by a manifold of unidirectional flux, which we call a ``dividing manifold'', and transition can be considered to correspond to crossing the dividing manifold. The stable (or forwards contracting) manifold of the NHS separates the extended state space into those initial conditions that will lead to a transition from those that will not. In this article, we rely on the stability and smoothness of NHS and their forward and backward contracting invariant submanifolds, developed by Fenichel in the early 1970s~\cite{fenichel1974asymptotic,fenichel1977asymptotic}.

In the one degree of freedom (1DoF) context, the NHS consists of a single trajectory. It is hyperbolic, generalising the hyperbolic equilibrium point at the maximum of the potential of the autonomous case. By an adaptation of~\cite{bishnani2003safety}, we describe how to compute this hyperbolic trajectory, and also its stable manifold. In the two degrees of freedom (2DoF) context, the NHS consists of a 2D set of trajectories (thus a 3D subset of the extended state space). It generalises from the autonomous Hamiltonian case the index-1 saddle points of the potential and the trajectories that oscillate about them along the ridge of the potential.

To expand on the previous paragraphs, a ``dividing manifold'' is constructed in the extended state-space; the ship can be deemed to capsize if and only if the trajectory crosses the dividing manifold in the appropriate direction.  It corresponds roughly to the boundary of the potential well for the upright state but has the property of being free of local recrossings (also referred to as ``no sneaky returns'' in~\cite{mackay1990flux}). It spans a NHS consisting of states that hover on the brink of capsizing, which we call a ``saddle manifold'', as it generalises the concept of a saddle point.  
Even though ``capsize rate'' is not a practical concept in the context of ship dynamics, we propose that identifying the boundary of a ship's safe operational states can improve existing active control and manoeuvring approaches.

First, we address a 1DoF model, namely passage over a potential barrier called the Eckart barrier borrowed from theoretical chemistry. Then, we consider a 2DoF model of roll and heave for a ship because these are considered to be the most significant degrees of freedom for studying the ship's motion near capsize \cite{thompson1996suppression}. 
The importance of nonlinear coupling of roll to other degrees of freedom for a ship such as pitch and heave was shown in~\cite{nayfeh1973nonlinear}, following up on work by~\cite{kinney1961unstable,paulling1959unstable} and dates back to Froude~\cite{froude1863remarks}. The pitch-roll model of~\cite{nayfeh1973nonlinear} explained the effect of pitch-to-roll frequency ratio on the saturation of pitch response and energy transfer to the roll DoF, and the resulting undesirable roll response.
In particular, we examine a roll-heave model of~\cite{thompson1996suppression}, present numerical algorithms to determine the associated saddle manifold and forward contracting manifold, and demonstrate how they can be used to classify initial conditions as safe or leading to capsize.

The structure of this paper is as follows: in Section~\ref{sec: theory and methods} the algorithms used to construct the normally hyperbolic submanifolds are presented. We begin by constructing the hyperbolic trajectories as a solution to a boundary value problem. In the case of a gap in the linearised attraction rates the associated centre manifold to this trajectory is then used as the NHS for the system. We also present a method to compute the stable manifolds for these objects. These methods are then applied in Section~\ref{sec: application to ship capsize} to a one degree of freedom model representing transition over a local potential maximum, and a two degree of freedom model representing escape from a potential well. In Section~\ref{sec: discussion and conclusions} the results of these methods on these systems are presented and the effectiveness of these methods are discussed.

\section{Theory and Methods}
\label{sec: theory and methods}
In this section we will consider a dynamical system defined by a system of ordinary differential equations
\begin{equation}
    \label{eq: ivp}
    \begin{cases}
            &\dot{\mathbf{x}} = \mathbf{f}(\mathbf{x},t), \\
            &\mathbf{x}(0) = \mathbf{x}_0 \in \mathbb R^d,
    \end{cases}
\end{equation}
where $\mathbf{f}(\mathbf{x},t)$  is a non-autonomous perturbation of an autonomous system, for example, of the form $\mathbf{f}_0(\mathbf{x}) + \mathbf{F}(t)$.  
For notational simplicity, we shall restrict attention to $\mathbf{f}$ of this form, but in general $\mathbf{F}$ could depend on $\mathbf{x}$ too.
In our applications in Section~\ref{sec: application to ship capsize}, $\mathbf{f}_0$ will be given by the evolution of a Hamiltonian system with an additional damping term; however, the methods outlined here are applicable to more general systems.

In the study of autonomous systems,  the stable and unstable manifolds to equilibria and periodic orbits serve as the skeleton from which the qualitative action of the flow on the phase space can be understood. 
Our approach generalises this approach to non-autonomous perturbations of dissipative systems, by making an extension of the definition of hyperbolic equilibria to a non-autonomous system.

In a first step, this leads to the idea of a \textit{hyperbolic trajectory} \cite{kuznetsov2012hyperbolic, madrid2009distinguished}. In general, these can not be expressed in a closed form, so we present a numerical algorithm to approximate such trajectories.

As a second step, it leads to the use of normally hyperbolic submanifolds for non-autonomous systems. We simplify here by considering only those normally hyperbolic submanifolds that arise as ``centre manifolds'' for hyperbolic trajectories. In a loose description, centre manifolds are invariant submanifolds whose trajectories remain relatively close to a hyperbolic trajectory in both directions in time, compared to others that diverge faster in forward or backward time. These centre manifolds will act as the saddle manifolds used to classify initial conditions as safe and unsafe. 
In the absence of dissipation and forcing, the centre manifold to a saddle point in a 2DoF Hamiltonian system is topologically a disk \cite{mackay1990flux}. In the case of dissipative and temporal perturbation, we would like to compute the centre manifold directly.


\subsection{Hyperbolic trajectories}
\label{sec: hyp traj 2dof}
A hyperbolic equilibrium of an autonomous system $\dot{\mathbf{x}}=\mathbf{f}(\mathbf{x})$, defined on a $n$-dimensional manifold $\mathbb{M}$ is a state $\mathbf{p} \in \mathbb{M}$ such that $\mathbf{f}(\mathbf{p})=0$ and the linearisation $D\mathbf{f}_{\mathbf{p}}$ at $\mathbf{p}$ has no eigenvalues on the imaginary axis.  Using local coordinates in $\mathbb{R}^n$,
the eigenvalue condition is equivalent to the linear map $L: C^1(\mathbb{R};\mathbb{R}^n)\to C^0(\mathbb{R};\mathbb{R}^n)$, defined by
\begin{equation}
    (\mathbf{L} \mathbf{\xi})(t) = \dot{\mathbf{\xi}}(t)-D\mathbf{f}_{\mathbf{p}} \mathbf{\xi}(t),
\end{equation} 
being invertible, where $C^k(\mathbb{R};\mathbb{R}^n)$ denotes the function space of $k$-times continuously differentiable functions $\mathbf{x}: \mathbb{R} \to \mathbb{R}^n$ with the $C^k$ norm.  Following \cite{bishnani2003safety},
we extend this idea to non-autonomous systems $\dot{\mathbf{x}}=\mathbf{f}(\mathbf{x},t)$ by defining a trajectory $\mathbf{x}:\mathbb{R} \to \mathbb{M}$ to be {\em hyperbolic} if $\mathbf{L}$ is invertible, where now $D\mathbf{f}_{\mathbf{p}}$ is replaced by $D\mathbf{f}_{(\mathbf{x}(t),t)}$, the matrix of $\mathbf{x}$-derivatives at $(\mathbf{x}(t),t)$.  This formulation requires the trajectory to remain in one local coordinate patch, which suffices for the present paper.  So, although the definition can be extended to general manifolds, we take $\mathbb{M}=\mathbb{R}^n$.

It is straightforward from the implicit function theorem that hyperbolic trajectories persist under $C^1$-small perturbation of a system 
(for a bounded invertible operator $\mathbf{L}$ between Banach spaces it is automatic that $\mathbf{L}^{-1}$ is bounded).
In the extended phase space $\mathbb{R}^d \times \mathbb{R}$ where time is added as an additional coordinate, a hyperbolic equilibrium $\mathbf{p}$ of the autonomous system corresponds to a hyperbolic trajectory of the extended system, which is given by
\begin{equation}
    \label{eq: extended ode}
    \begin{pmatrix}
    \dot{\mathbf{x}} \\
    \dot{t}
\end{pmatrix} = \begin{pmatrix}
    \mathbf{f}(\mathbf{x},t) \\
    1
\end{pmatrix}.
\end{equation}
This  hyperbolic trajectory  persists for sufficiently small forcing values, measured in the $C^1$ norm, to a locally unique nearby hyperbolic trajectory for the forced system~\cite{bishnani2003safety}.  
A key property of the hyperbolic trajectory is that it is a bounded function of time, in both directions of time. Indeed, it is the unique trajectory that remains forever in a uniform neighbourhood.  We use this property to {help} find the hyperbolic trajectory. 

If the forcing function is specified on only a bounded time-interval then the hyperbolic trajectory is not uniquely determined, but assuming the forcing function is bounded then the ambiguity decays exponentially from the ends of the interval into its interior, so the hyperbolic trajectory is accurately determined in a subinterval excluding neighbourhoods of the ends by choosing any bounded extension of the forcing, e.g.~zero.

Throughout this section, we will take $\mathbf{p}$ to be a hyperbolic equilibrium of $\mathbf{f}_0$, and let $\{\mathbf{u}_{i}\}_{i=1,\dotsc, N_-}$ and $\{\mathbf{v}_{i}\}_{i=1,\dotsc, N_+}$ be bases for the forwards and backwards contracting subspaces respectively, where $N_-,N_+$  are the dimensions of the forwards/backwards contracting subspaces of $D\mathbf{f}_0(\mathbf{p})$. In the case of distinct real eigenvalues one could choose eigenvectors corresponding to the negative and positive eigenvalues of $D\mathbf{f}_0(\mathbf{p})$ respectively.
Let $\mathbf{P}$ be the change of basis matrix given by
\begin{equation}
    \mathbf{P}_{i,j} = \begin{cases}
        (\mathbf{u}_{j})_i, \quad &\text{if } 1\le j \le N_-, \\
        (\mathbf{v}_{j- N_-})_i, \quad &\text{if } N_-< j \le d.
    \end{cases}
\end{equation}
Then given a solution of the autonomous version of Eqn.~\eqref{eq: extended ode} which is initially close to $x_0$, the displacement of the solution from $\mathbf{p}$ would expand (or contract) along the corresponding subspaces. If the forcing acts as a small 
$C^1$-perturbation of this system, the exponential growth of the autonomous system would dominate over a large time interval. These directions therefore act as approximate  forwards contracting/backwards contracting directions for the non-autonomous system, and we propose that for sufficiently small forcing, the hyperbolic trajectory $\mathbf{x}(\cdot)$ approximately solves the boundary value problem for large $T>0$:
\begin{equation}
    \label{eq: hyperbolic BVP}
        \begin{cases}
            &\dot{\mathbf{x}} = \mathbf{f}(\mathbf{x},t), \\
            &(\mathbf{P}^{-1}(\mathbf{x}(-T) -\mathbf{p}))_j = 0 , \quad \text{for } j = 1,\dotsc,N_-, \\
            &(\mathbf{P}^{-1}(\mathbf{x}(T) -\mathbf{p}))_j = 0,  \quad \text{for } j = N_- + 1,\dotsc,d.
    \end{cases}
\end{equation}
The first boundary condition specifies that the solution does not grow too large along the unstable directions of the autonomous system over a large time period, while the second specifies that the solution does not grow too large along the stable directions as time decreases. Both forward and backward directions are needed to set up this problem, as there are many solutions to Eqn.~\eqref{eq: extended ode} which remain bounded as time increases or as time decreases; the hyperbolic trajectory is the only solution which does both.

We now demonstrate how to compute the hyperbolic trajectory analytically and outline the numerical method, in the case where there are both positive and negative eigenvalues.
\subsubsection{Hyperbolic trajectory of a linear system}
\label{sec: linear hyperbolic trajectory}
It is sometimes possible to express the hyperbolic trajectory of a linear system of ODEs analytically, and this serves as both a useful exercise in understanding the hyperbolic trajectories and also a step in the algorithm described in Section \ref{sec: algorithm hyp traj}. Letting $\mathbf{A} = D\mathbf{f}_0(\mathbf{p})$, the \textit{linearised system} is defined to be
\begin{equation}
    \label{eq: linearised system}
    \begin{cases}
    \dot{\mathbf{y}} = \mathbf{A}\mathbf{y} + \mathbf{F}(t), \\ \quad \mathbf{y}(0) = \mathbf{y}_0.
    \end{cases}
\end{equation}
When a closed form general solution to \eqref{eq: linearised system} can be found, the hyperbolic trajectory is given by choosing the initial condition $\mathbf{y}_0$ so that $\mathbf{y}(t)$ remains bounded in both forwards and backwards time. 
We present the case for $d = 2$, $\mathbf{A} = \begin{pmatrix}
    -\mu_- & 0 \\
    0 & \mu_+
\end{pmatrix}$ and $\mathbf{F}(t) = ( -a\cos(\omega t), a\cos (\omega t))^T$, which is relevant to the application to the model in Section~\ref{subsec: one dof model}.
In this case it follows that the general solution to \eqref{eq: linearised system} is given by
\begin{equation}
\mathbf{y}(t) = \begin{pmatrix}
    e^{-\mu_- t} & 0 \\
    0 & e^{\mu_+ t}
\end{pmatrix}\mathbf{y}_0 + \begin{pmatrix}
    \int_0^t -ae^{-\mu_- s} \cos( \omega s) ds \\
    \int_0^t ae^{\mu_+ s}\cos(\omega s) ds 
\end{pmatrix}.
\end{equation}
The integral term can be simplified using that:
\begin{equation}
    \int_0^t e^{\mu s} \cos (\omega s) \; ds = \left(\frac{1}{1 + \frac{\mu^2}{\omega^2}} \right) \left[ e^{\mu t}\left( \frac{1}{\omega} \sin (\omega t) + \frac{\mu}{\omega^2} \cos (\omega t)\right) - \frac{\mu}{\omega^2} \right] .
\end{equation}
Letting $\alpha_{\pm} = \dfrac{a}{1 + \frac{\mu_{\pm}^2}{\omega^2}}$, the general solution to can be written as
\begin{equation}
    \label{eq: expanded general solution}
    \mathbf{y}(t) =   \begin{pmatrix}
        e^{-\mu_- t} (\mathbf{y}_0)_1 \\
        e^{\mu_+ t} (\mathbf{y}_0)_2
    \end{pmatrix}  + \begin{pmatrix}
     -a\alpha_- \left[\frac{1}{\omega} \sin (\omega t) + \frac{\mu_-}{\omega^2} \cos (\omega t) \right]  +a e^{-\mu_-t}   \frac{\mu_-}{\omega^2} \alpha_- \\
    a\alpha_+\left[\frac{1}{\omega} \sin (\omega t) + \frac{\mu_+}{\omega^2} \cos (\omega t) \right]  -a e^{-\mu_+ t} \frac{\mu_+}{\omega^2}\alpha_+
    \end{pmatrix}.
\end{equation}
By choosing the initial condition so that all terms involving an exponential vanish, we obtain the hyperbolic trajectory:
\begin{equation}
    \mathbf{y}_{hyp}(t) =       \begin{pmatrix}
     -a \alpha_-\left[\frac{1}{\omega} \sin (\omega t) + \frac{\mu_-}{\omega^2} \cos (\omega t) \right]  \\
    a \alpha_+ \left[\frac{1}{\omega} \sin (\omega t) + \frac{\mu_+}{\omega^2} \cos (\omega t) \right]
    \end{pmatrix} .
\end{equation}
In general we will not be able to find a closed form expression for the hyperbolic trajectory of a nonlinear system, and hence a numerical algorithm is necessary. We propose an iterative algorithm, which requires a good initial guess for the hyperbolic trajectory. For this we may use the hyperbolic trajectory to the corresponding linearised system about the saddle point. 
Alternatively, we may use the constant trajectory $\mathbf{x}(t)=\mathbf{p}$ as an initial guess.
\subsubsection{Algorithm for finding hyperbolic trajectory}
\label{sec: algorithm hyp traj}
Given a forcing function specified over some time-interval $I$ we present an algorithm to numerically determine the corresponding hyperbolic trajectories of the system, using the Newton-Raphson method (subject to suitable boundary conditions). The procedure to find the hyperbolic trajectory for forcing of a saddle point $\mathbf{p}$ goes as follows:
\begin{enumerate}
    \item Without loss of generality  we assume that $I = [-T,T]$ and let $N$ be the number of desired points of the hyperbolic trajectory. It is best that $N$ is large and we shall assume that it is greater than 3. Let $\mathbf{t} \in \mathbb{R}^N$ be a discretisation of the time interval, that is $-T = t_1 < t_2 < \dotsb < t_{N-1} < t_N = T$. We will assume that this discretisation is equally spaced, so $t_{i+1} = t_i + \delta t$ with $\delta t = \frac{2T}{N}$.
    \item Let $\mathbf{x}: I \to \mathbb{R}^d$ be an initial guess for the hyperbolic trajectory. This could be the corresponding hyperbolic trajectory for the linearised system with the given forcing function; an example of this was given in Section \ref{sec: linear hyperbolic trajectory}. It is also possible to choose the saddle point, $\mathbf{x} : t \mapsto \mathbf{p}$.
    \item Let $\mathbf{X} \in (\mathbb{R}^d)^N$ be given by $\mathbf{X}_i = \mathbf{x}(t_i) - \mathbf{p}$, for $i=1,\dotsc, N$. This is the displacement of the hyperbolic trajectory from the saddle point $\mathbf{p}$. We define $\mathbf{\Phi}: (\mathbb{R}^d)^N \to (\mathbb{R}^d)^N$ by
    \begin{equation} 
    \mathbf{\Phi}(\mathbf{Y})_i = \begin{cases}
        & \phi_{\delta t}(\mathbf{Y}_{i}+ \mathbf{p},t_{i}) - \mathbf{Y}_{i+1} - \mathbf{p},\quad \text{ for } i = 1,\dotsc, N-1, \\
        &\mathbf{B}(\mathbf{Y}_1,\mathbf{Y}_N), \quad  \text{ for } i =N,
    \end{cases}
    \end{equation}
    where $\phi$ is the flow map for the non-linear differential equation with forcing, and the boundary condition $\mathbf{B}(\mathbf{Z}_1,\mathbf{Z}_2)$ is given by:
    \begin{equation}
        \mathbf{B}(\mathbf{Z}_1,\mathbf{Z}_2)_j = \begin{cases}
            &(\mathbf{P}^{-1}(\mathbf{Z}_1 -\mathbf{p}))_j,   \quad \text{for } j = 1,\dotsc,N_-, \\
            &(\mathbf{P}^{-1}(\mathbf{Z}_2 -\mathbf{p}))_j,   \quad \text{for } j = N_- + 1,\dotsc,d,
        \end{cases}
    \end{equation}
    \item A vector $\mathbf{Y} \in (\mathbb{R}^d)^N$ with $\mathbf{\Phi}(\mathbf{Y}) = \mathbf{0}$ will correspond to a discretised hyperbolic trajectory given by $(\mathbf{Y}_i + \mathbf{p})_{i=1,\dotsc N}$. We approximate the derivative of the map $\phi_{\delta t }(\cdot + \mathbf{p})$ by the derivative of the flow map of the corresponding linear system, which is given by 
    $e^{\delta t \mathbf{A}}$, where $\mathbf{A} = D\mathbf{f}(\mathbf{p})$. Therefore beginning with $\mathbf{X}^0 = \mathbf{X}$ we apply the iterative scheme
    \begin{equation} \mathbf{X}^{k+1} = \mathbf{X}^k - \mathbf{M}^{-1}\mathbf{\Phi}(\mathbf{X}^k) ,\end{equation}
    where
    \begin{equation}
    \label{eq: approx DGS derivative 1}
    M= \begin{pmatrix}
        e^{\delta t \mathbf{A}} & -\mathbf{I}             & \mathbf{0}                & \cdots  & \mathbf{0}& \cdots  & \mathbf{0}      & \mathbf{0} \\
        \mathbf{0}              & e^{\delta t \mathbf{A}} & -\mathbf{I}                & \cdots & \mathbf{0}    & \cdots  & \mathbf{0}  & \mathbf{0} \\
        \mathbf{0}              & \mathbf{0}              & e^{\delta t \mathbf{A}}    & \cdots & \mathbf{0}    & \cdots  & \mathbf{0}  & \mathbf{0} \\
        \vdots         & \vdots         &  \vdots           & \vdots & \ddots & \cdots  & \vdots & \vdots \\
        \mathbf{B}_2          & \mathbf{0}              & \mathbf{0}                 & \mathbf{0} & \mathbf{0} & \cdots  & \mathbf{0} & \mathbf{B}_1
    \end{pmatrix} ,
\end{equation}
and $\mathbf{B}_1, \mathbf{B}_2 \in \mathbb R ^{N_+\times d} $ are given by:
\begin{equation}
    (\mathbf{B}_1)_{i,j} = \begin{cases}
        \mathbf({P}^{-1})_{j,i},  \quad \text{for } i =1.\dotsc N_-,j = 1,\dotsc d, \\
        0,  \quad \text{otherwise},
    \end{cases}
\end{equation}
\begin{equation}
    (\mathbf{B}_2)_{i,j} = \begin{cases}
        \mathbf({P}^{-1})_{j,i},  \quad \text{for } i =N_-+1,\dotsc d, , j = 1,\dotsc d, \\
        0,  \quad \text{otherwise}. 
    \end{cases} 
\end{equation}
\item This iteration procedure is continued until it is deemed to have converged using the absolute convergence criteria of $\left\vert\left\vert \mathbf{X}^{k} - \mathbf{X}^{k-1} \right\vert\right\vert < \varepsilon_c$ and $\left\vert\left\vert \mathbf{\Phi}(\mathbf{X}^k) \right\vert\right\vert < \varepsilon_F$ from \cite[pp.~302]{parker2012practical}. In practice we do not compute $\mathbf{M}^{-1}$, instead a LUP decomposition of $\mathbf{M}$ is computed once and stored, and then we solve the equation $\mathbf{M}(\delta \mathbf{X}^k) = -\mathbf{\Phi}(\mathbf{X}^k)$ and iterate $\mathbf{X}^{k+1} =\mathbf{X}^{k} +  \delta \mathbf{X}^k$.
\end{enumerate}

This method is in principle extendible to forcing that depends on the spatial variable $\mathbf{x}$. In this case the initial guess for the hyperbolic trajectory given by the linearised system is no longer valid, and the spatial derivative of $\phi_{\delta t} (\cdot + \mathbf{p})$ may no longer be well approximated by $e^{\delta t \mathbf{A}}$.
To overcome these difficulties one could approximate the forcing $\mathbf{F}(\mathbf{x},t)$ by $D_{\mathbf{x}}\mathbf{F}(\mathbf{p},t) \mathbf{x} + \mathbf{F}(\mathbf{p},t)$. Letting $\mathbf{B}(t) = D_{\mathbf{x}}\mathbf{F}(\mathbf{p},t)$, the linearised system then becomes:
\begin{equation}\dot{\mathbf{y}} = (\mathbf{A}+\mathbf{B}(t))\mathbf{y} + \mathbf{F}(\mathbf{p},t). \end{equation}
From here it is now possible to approximate the derivative of $\phi_{\delta t}(\cdot + \mathbf{p})$ by the derivative of the general solution to linearised solution.
Furthermore, if a general solution is available to the linearised system then this can be used to obtain an initial guess for the hyperbolic trajectory.
For the remainder of this section we will again concern ourselves only with forcing of the form $\mathbf{F}(t)$.

\subsection{Centre, stable and unstable manifolds}
A hyperbolic trajectory has stable and unstable manifolds, defined as the sets of points whose forwards, respectively backwards, trajectory converges to it.  They form injectively immersed smooth submanifolds of the extended state-space.  Furthermore, if there is one direction of forward contraction that is stronger than the rest, the stable manifold contains a strong stable manifold of dimension one and a complementary submanifold of codimension-one that is often called a centre manifold. 

The existence of stable and centre manifolds for the hyperbolic trajectory of the non-linear system can be seen from \cite{jolly2005computation,mancho2003computation}. These manifolds are of interest as they  persist under weak dissipation and external forcing. In the context of passage over a saddle point, the centre manifold serves as a \textit{saddle manifold} \cite{bujorianu2021new}, with codimension 1 forward and backward contracting manifolds.  This saddle manifold can be spanned by a \textit{dividing manifold}, a codimension 1 submanifold that is transverse to the vector field, except on the saddle manifold. 
Note that  there is not a unique choice for the dividing manifold, as it is defined by transversality rather than invariance, but the resulting ambiguity is only in the precise time at which one declares the transition to occur. The forward contracting manifold serves to partition the extended state space into the regions whose trajectories pass through the dividing manifold and so over the centre manifold, and those which do not.

Given a forcing function known over some time interval $I$ we now present a modification to the algorithm in Section \ref{sec: algorithm hyp traj} to numerically determine centre and stable manifolds corresponding to the saddle point $\mathbf{p}$.
\subsubsection{Centre manifold algorithm}
\label{sec: algorithm center manifold}
Let $\mathbf{A}= \mathbf{Df}_0(\mathbf{p})$ and $\{\mathbf{u}_{i}\}_{i=1,\dotsc, N_-}$, $\{\mathbf{w}_i\}_{i = 1,\dotsc N_C} $ and $\{\mathbf{v}_{i}\}_{i=1,\dotsc, N_+}$ denote real basis vectors for the stable/centre/unstable directions respectively. In the context of the example in Section~\ref{subsec: two dof model}, the centre eigenvalues might not be purely imaginary, however they will not expand/contract as fast as the stable and unstable eigenvalues. In the following we assume that all eigenvalues along the centre subspace have negative real part, which is the case in Section \ref{subsec: two dof model} for a Hamiltonian system with weak dissipation.
Let $\mathbf{P}$ be the change of basis matrix given by:
\begin{equation}
    \mathbf{P}_{i,j} = \begin{cases}
        (\mathbf{u}_{j})_i, \quad &\text{if } 1\le j \le N_-, \\
        (\mathbf{w}_{j-N_-})_i, \quad &\text{if } N_-< j \le N_- +N_C, \\
        (\mathbf{v}_{j-N_C- N_-})_i, \quad &\text{if } N_- + N_C < j \le d.
    \end{cases}
\end{equation}
As in Section \ref{sec: algorithm hyp traj}, first the linearised system Equation (\ref{eq: linearised system}) is examined. Letting $\mathbf{x}_{hyp,lin}$  denote the hyperbolic trajectory of the linearised system, it is clear that the centre manifold is given by 
\begin{equation}
    \label{eq: linearised centre manifold}
    \left\{\left(\mathbf{x}_{hyp,lin}(t)+ \sum_{i=1}^{N_c} e^{t \mathbf{A}}(q_i \mathbf{w}_i),t \right) :\: t,q_i \in \mathbb{R}, \; i = 1,\dotsc N_C  \right\}.
\end{equation}
Letting $\mathbf{y}_C(t;\mathbf{q}) = \sum_{i=1}^{N_c} e^{t \mathbf{A}}(q_i \mathbf{w}_i) $ for $\mathbf{q} \in \mathbf{R}^{N_C}$, the above is equivalent to $\mathbf{x}_{hyp,lin}(t) + \mathbf{y}_C(t;\mathbf{q})$ being a trajectory along the centre manifold.  As the centre manifold to the linearised system is tangent to the centre manifold of the non-linear system in the case of no external forcing, the centre manifold to the non-linear system $\dot{\mathbf{x}} = \mathbf{f}_0(\mathbf{x})$ is  locally given  by a graph over the centre subspace of $\mathbf{A}$. The centre manifold depends smoothly upon the vector field, and so for small enough forcing, the centre manifold to the non-linear forced system will also be given locally by a graph over the centre subspace. In particular, the orthogonal projection of the centre manifold onto the affine subspace of $\mathbf{A}$ through $\mathbf{p}$ is \emph{bijective} in a neighbourhood of the hyperbolic trajectory for sufficiently small forcing values.

The algorithm to find the centre manifold for the non-linear and forced system goes as follows:
\begin{enumerate}
    \item Let $I$, $N$, $\mathbf{t}$, $\delta t$ be as in Step 1 of Section \ref{sec: algorithm hyp traj}, and let $\mathbf{X}_{hyp}$ be the hyperbolic trajectory computed according to Section \ref{sec: algorithm hyp traj}. Let $J \in \{1,\dotsc, N \}$ be an index which is not too close to either end, and $\mathbf{q} \in \mathbb{R}^{N_C}$. 
    \item Let $\mathbf{Y} \in (\mathbb{R}^d)^N$ be given by $\mathbf{Y}_i =  \mathbf{y}_C(t;\mathbf{q})(\mathbf{t}_i)$, for $i=1,\dotsc, N$. This is the displacement of a centre manifold trajectory from the hyperbolic trajectory. We define $\mathbf{\Phi_C}: (\mathbb{R}^d)^N \to (\mathbb{R}^d)^N$ by
    \begin{equation} 
    \mathbf{\Phi_C}(\mathbf{Y})_i = \begin{cases}
        &  \phi_{\delta t}(\mathbf{Y}_{i}+\mathbf{X}_{hyp,i} ,t_{i}) - \mathbf{Y}_{i+1} - \mathbf{X}_{hyp,i+1}, \quad \text{ for } i = 1,\dotsc, N-1, \\
        &\mathbf{B_C}(\mathbf{Y}_1,\mathbf{Y}_J,\mathbf{Y}_N), \quad  \text{ for } i =N.
    \end{cases}
    \end{equation}
    where $\phi$ is the flow map for the non-linear differential equation with forcing, and the boundary condition $\mathbf{B_C}(\mathbf{Z}_1,\mathbf{Z}_2,\mathbf{Z}_3)$ is given by:
    \begin{equation}
        \mathbf{B_C}(\mathbf{Z}_1,\mathbf{Z}_2,\mathbf{Z}_3)_j = \begin{cases}
            &(\mathbf{P}^{-1}(\mathbf{Z}_1 -(\mathbf{X}_{hyp})_1))_j ,  \quad \text{for } j = 1,\dotsc,N_-, \\
            &(\mathbf{P}^{-1}(\mathbf{Z}_2 -(\mathbf{X}_{hyp})_J))_j - q_{j-N_-}, \text{for } j - N_- = 1 ,\dotsc,N_C, \\
            &(\mathbf{P}^{-1}(\mathbf{Z}_3 -(\mathbf{X}_{hyp})_N))_j  , \quad \text{for } j = N_- +N_C + 1,\dotsc,d.
        \end{cases}
    \end{equation}
The boundary condition $\mathbf{B_C} = \mathbf{0}$ can be interpreted as ensuring that the solution does not blow up along the forwards and backwards contracting directions, and along the centre directions it has a displacement of $\sum_{i=1}^{N_C} q_i \mathbf{w}_i$ from the hyperbolic trajectory, when projected into the centre subspace.
    \item A discretised trajectory along the centre manifold will correspond to a $\mathbf{Y}+\mathbf{X}_{hyp}$ where $\mathbf{Y} \in (\mathbb{R}^d)^N$ with $\mathbf{\Phi_C}(\mathbf{Y}) = 0$; the remaining steps continue as in Steps 3 and 4 of the algorithm in Section \ref{sec: algorithm hyp traj}, where now the derivative of $\mathbf{\Phi_C}$ is approximated by:
    \begin{equation}
    \label{eq: approx DFC derivative 1}
    \begin{pmatrix}
        e^{\delta t \mathbf{A}} & -\mathbf{I}             & \mathbf{0}                & \cdots  & \mathbf{0}& \cdots  & \mathbf{0}      & \mathbf{0} \\
        \mathbf{0}              & e^{\delta t \mathbf{A}} & -\mathbf{I}                & \cdots & \mathbf{0}    & \cdots  & \mathbf{0}  & \mathbf{0} \\
        \mathbf{0}              & \mathbf{0}              & e^{\delta t \mathbf{A}}    & \cdots & \mathbf{0}    & \cdots  & \mathbf{0}  & \mathbf{0} \\
        \vdots         & \vdots         &  \vdots           & \vdots & \ddots & \cdots  & \vdots & \vdots \\
        \mathbf{B}_1          & \mathbf{0}              & \mathbf{0}                 & \mathbf{B_J} & \mathbf{0} & \cdots  & \mathbf{0} & \mathbf{B}_N
    \end{pmatrix},
\end{equation}
where $\mathbf{B_J}$ is a block matrix covering the $dJ$th to the $d(J+1)$th columns, and
$\mathbf{B}_1, \mathbf{B}_J,\mathbf{B}_N \in \mathbb R ^{d \times d} $ are given by:
\begin{align}
    (\mathbf{B}_1)_{i,j} &= \begin{cases}
        \mathbf({P}^{-1})_{j,i},  \quad \text{for } i =1,\dotsc N_-,, j = 1,\dotsc d, \\
        0,  \quad \text{otherwise} ,
    \end{cases}\\
    (\mathbf{B}_J)_{i,j} &= \begin{cases}
        \mathbf({P}^{-1})_{j,i},  \quad \text{for } i =N_-+1,\dotsc N_- + N_C, , j = 1,\dotsc d, \\
        0,  \quad \text{otherwise} ,
    \end{cases}\\
    (\mathbf{B}_N)_{i,j} &= \begin{cases}
        \mathbf({P}^{-1})_{j,i},  \quad \text{for } i =N_-+N_C + 1,\dotsc d, , j = 1,\dotsc d, \\
        0,  \quad \text{otherwise} .
    \end{cases}
\end{align}
\end{enumerate}

\subsubsection{Stable and unstable manifold algorithm}
\label{subsec: algorithm stable manifold}
The approach to find the centre manifold in Section \ref{sec: algorithm center manifold} can be further adapted to find the stable manifold. For this, we observe that for the linearised system the stable manifold to the centre manifold is given by:
\begin{equation}
    \label{eq: linearised stable manifold}
    \begin{split}
    \left\{\left(\mathbf{x}_{hyp,lin}(t) + \sum_{i=1}^{N_c} e^{t \mathbf{A}}(q_i \mathbf{w}_i) + \sum_{j=1}^{N_-} e^{t \mathbf{A}}(p_j \mathbf{u}_j),t\right)\right. \\ 
    \left.:\: t,q_i,p_j \in \mathbb{R}, \; i = 1,\dotsc N_C, j =1 ,\dotsc N_-  \right\}.    
    \end{split}
\end{equation}
As for the centre manifold, the stable manifold to the non-linear autonomous system $\dot{\mathbf{x}} =  \mathbf{f}_0(\mathbf{x})$ is locally given by a graph over the affine subspace $\mathbf{p} + \mathrm{Span}(\mathbf{u}_1,\dotsc,\mathbf{u}_{N_-}, \mathbf{w}_1,\dotsc \mathbf{w}_{N_C})$, and as the stable manifold depends smoothly upon the vector field, for each $t_0 \in \mathbb{R}$ it is locally given by a graph over the affine subspace $\mathbf{x}_{hyp}(t_0) + \mathrm{Span}(\mathbf{u}_1,\dotsc,\mathbf{u}_{N_-}, \mathbf{w}_1,\dotsc \mathbf{w}_{N_C})$.

Letting $\mathbf{y}_S(t,\mathbf{q}) := \sum_{i=1}^{N_c} e^{t \mathbf{A}}(q_i \mathbf{w}_i) + \sum_{j=1}^{N_-} e^{t \mathbf{A}}(q_{j+N_C} \mathbf{u}_j) $, for $\mathbf{q} \in \mathbb{R}^{N_- + N_C}$, the algorithm to find the stable manifold is given by:
\begin{enumerate}
    \item Let $I$, $N$, $\mathbf{t}$, $\delta t$ be as in Step 1 of Section \ref{sec: algorithm hyp traj}, and let $\mathbf{X}_{hyp}$ be the hyperbolic trajectory computed according to Section \ref{sec: algorithm hyp traj}. Let $J \in \{1,\dotsc, N \}$ be an index which is not too close to either end, and $\mathbf{q} \in \mathbb{R}^{N_C+N_-}$. 
    \item Let $\mathbf{Y} \in (\mathbb{R}^d)^N$ be given by $\mathbf{Y}_i =\mathbf{y}_S(t;\mathbf{q})(\mathbf{t}_i)$ for $i=1,\dotsc, N$. We define $\mathbf{\Phi_S}: (\mathbb{R}^d)^N \to (\mathbb{R}^d)^N$, by
    \begin{align} 
    &\mathbf{\Phi_S}(\mathbf{Y})_i = \notag \\
    &\begin{cases}
        &\phi_{\delta t}(\mathbf{Y}_{i}+\mathbf{X}_{hyp,i} ,t_{i}) - \mathbf{Y}_{i+1} - \mathbf{X}_{hyp,i+1}, \quad \text{for} \; i = 1,\ldots, N-1, \\
        &\mathbf{B_S}(\mathbf{Y}_J,\mathbf{Y}_N) , \quad \text{for} \; i =N,
    \end{cases}
    \end{align}
    where $\phi$ is the flow map for the non-linear differential equation with forcing, and the boundary condition $\mathbf{B_S}(\mathbf{Z}_1,\mathbf{Z}_2)$, is given by:
    \begin{equation}
        \mathbf{B_S}(\mathbf{Z}_1,\mathbf{Z}_2)_j = \begin{cases}
            &(\mathbf{P}^{-1}(\mathbf{Z}_1 -\mathbf{X}_{hyp,J}))_ - q_{j}  ,\quad \text{for } j = 1 ,\dotsc,N_-+N_C, \\
            &(\mathbf{P}^{-1}(\mathbf{Z}_2 -\mathbf{X}_{hyp,_N}))_j , \quad \text{for } j = N_- +N_C + 1,\dotsc,d.
        \end{cases}
    \end{equation}
The boundary condition $\mathbf{B_S} = \mathbf{0}$, should be interpreted as ensuring that the solution does not blow up along the backwards contracting directions, and along the centre and stable directions it has a displacement of $\sum_{i=1}^{N_C} q_i \mathbf{w}_i + \sum_{j=1}^{N_-} q_{j+N_C} \mathbf{u}_j$ from the hyperbolic trajectory, when projected into the centre and stable subspaces.
    \item A  trajectory along the stable manifold will correspond to  $\mathbf{Y}+\mathbf{X}_{hyp}$ where $\mathbf{Y} \in (\mathbb{R}^d)^N$ with $\mathbf{\Phi_S}(\mathbf{Y}) = 0$ ; the remaining steps continue as in Steps 3 and 4 of the algorithms in Sections \ref{sec: algorithm hyp traj} and \ref{sec: algorithm center manifold}, where now the derivative of $\mathbf{\Phi_S}$ is approximated by:
    \begin{equation}
    \label{eq: approx DFC derivative 2}
    \begin{pmatrix}
        e^{\delta t \mathbf{A}} & -\mathbf{I}             & \mathbf{0}                & \cdots  & \mathbf{0}& \cdots  & \mathbf{0}      & \mathbf{0} \\
        \mathbf{0}              & e^{\delta t \mathbf{A}} & -\mathbf{I}                & \cdots & \mathbf{0}    & \cdots  & \mathbf{0}  & \mathbf{0} \\
        \mathbf{0}              & \mathbf{0}              & e^{\delta t \mathbf{A}}    & \cdots & \mathbf{0}    & \cdots  & \mathbf{0}  & \mathbf{0} \\
        \vdots         & \vdots         &  \vdots           & \vdots & \ddots & \cdots  & \vdots & \vdots \\
        \mathbf{0}         & \mathbf{0}              & \mathbf{0}                 & \mathbf{B_J} & \mathbf{0} & \cdots  & \mathbf{0} & \mathbf{B}_N
    \end{pmatrix},
\end{equation}
where $\mathbf{B_J}$ is a block matrix covering the $dJ$th to the $d(J+1)$th columns, and
$\mathbf{B}_J,\mathbf{B}_N \in \mathbb R ^{d \times d} $ are given by:
\begin{equation}
\begin{split}
    (\mathbf{B}_J)_{i,j} = \begin{cases}
        \mathbf({P}^{-1})_{j,i},  \quad \text{for } i =1,\dotsc N_- + N_C, , j = 1,\dotsc d, \\
        0,  \quad \text{otherwise} ,
    \end{cases}\\
    (\mathbf{B}_N)_{i,j} = \begin{cases}
        \mathbf({P}^{-1})_{j,i} , \quad \text{for } i =N_-+N_C + 1,\dotsc d, , j = 1,\dotsc d, \\
        0 , \quad \text{otherwise}.
    \end{cases}
    \end{split}
\end{equation}
\end{enumerate}
By making the appropriate changes to the boundary conditions and initial trajectory choice, we can also apply this method to construct the unstable manifold.
\subsection{Parameterising the stable and centre manifolds}
\label{sec: parameterising manifolds}
The method in Section \ref{subsec: algorithm stable manifold} provides a natural immersion of the stable manifold into the extended state space, which is given by the mapping $\iota_S :(\mathbf{q},t_0) \mapsto (\mathbf{y}(t_0),t_0)$ which maps $(\mathbf{q},t_0)$ to the trajectory $\mathbf{y}$ outputted by the algorithm (with inputs $(\mathbf{q},t_0)$) evaluated at time $t_0$. That is
$\iota_S(\mathbf{q},t_0) = \mathbf{Y}_J$, where $\mathbf{t}_J = t_0$. Letting $\Omega_S = \{(\mathbf{q},\tau): \mathbf{q} \in [-R,R]^{N_-+N_C}, \tau \in (\tau_{-}, \tau_{+}) \}$,
the immersion then maps $\iota_{S} : \Omega_S \to \mathbb{R}^d \times \mathbb{R}$. 
By solving the boundary value problem for a sample of points in $\Omega_S$, we can then reconstruct this immersion through interpolation, to obtain an approximate parametrisation  of the stable manifold defined on the whole of $\Omega_{S,\pm}$. 
A similar method using the algorithm in Section \ref{sec: algorithm center manifold} can also be used to construct an approximate parametrisation of the centre manifold $\iota_{C} : \Omega_{C} \to \mathbb{R}^d \times \mathbb{R}$ where $\Omega_{C} = \{(\mathbf{q},\tau): \mathbf{q} \in [-R,R]^{N_C}, \tau \in (\tau_{-}, \tau_{+})\}$.
 These interpolations can be used to artificially increase the number of sampled points from the centre manifold, and  for visualisations of these manifolds. 
 
In implementation, some modifications are made to the algorithms  in Sections  \ref{sec: algorithm center manifold} and \ref{subsec: algorithm stable manifold} to improve efficiency. As only one value of the trajectory is used, the BVP can be considered over a shorter time interval, $[\tau_{min}, \tau_{max}]$ instead of $[-T,T]$.  Specifically, we express this interval about the ``initial'' time $t_j$ by
$\tau_{min} = t_j - \Delta T_-,$ and $\tau_{max} = t_j + \Delta T_+$,
where $\Delta T_\pm$ are to be chosen appropriately. These must be chosen large enough so that one would expect the stable and unstable directions to see growth over their respective interval but small enough to reduce the possibility of blow up and to decrease computational cost.

One such choice is to choose these time intervals in relation to the strength of contraction/expansion of the system. Letting $\lambda_-$ be the smallest (most negative) eigenvalue in the stable subspace and $\lambda_+$ be the greatest (most positive) eigenvalue in the unstable subspace, $\Delta T_\pm$ were chosen such that:
\begin{equation}e^{- \lambda_+ \Delta T_-}  = e^{ \lambda_- \Delta T_+} = 10. \end{equation} 
\section{Applications}
\label{sec: application to ship capsize}
We now apply our methods to study the passage over a saddle point in two perturbed Hamiltonian systems.
\subsection{One DoF model}
\label{subsec: one dof model}
\subsubsection{Autonomous undamped model}
\label{subsec: 1dof autonomous undamped model}
In a Hamiltonian system with one degree of freedom the problem of escaping over a saddle point reduces to escaping over a local maximum of the energy potential. As such a system has a two-dimensional phase space it is easier to visualise than higher dimensional systems and so the central idea behind the method presented in this paper is first explored within this setting. 

\begin{figure}[hbtp]
    \centering
    \includegraphics[width=0.7\textwidth]{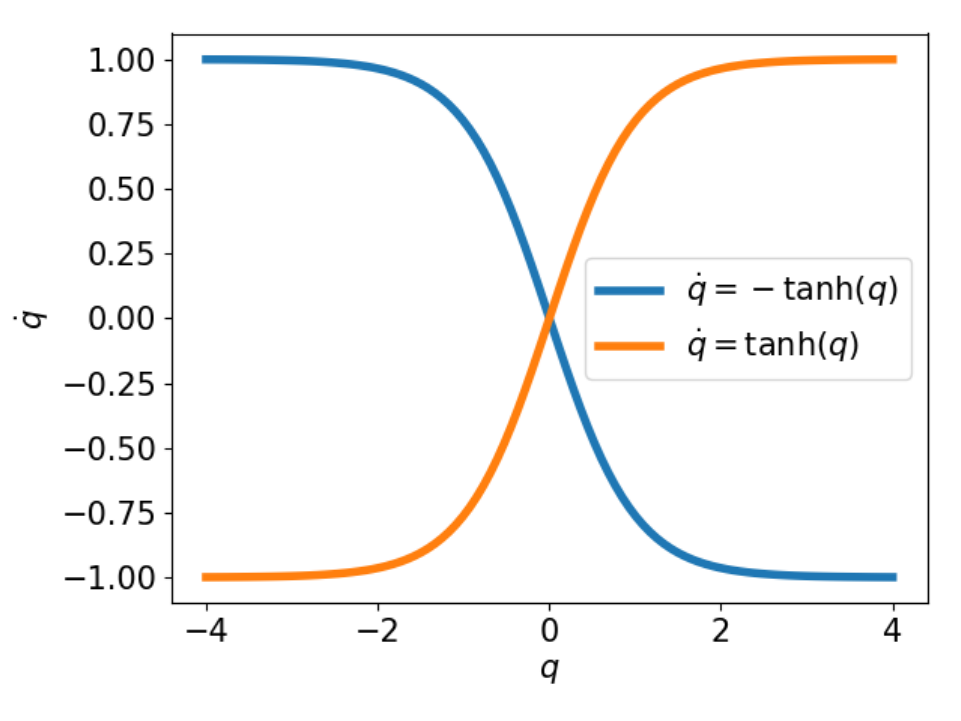}
    \caption{The stable and unstable manifolds to the origin for the  conservative system with $A = 1/2, I = 1$. The system is then $\ddot{q} = \tanh (q) \mathrm{sech}^2 (q)$. The stable manifold is given by the graph of $\dot{q} = -\tanh q$ and the unstable manifold by $\dot{q} =\tanh q$}
    \label{fig:stableManifoldConservative}
\end{figure}

We shall first study a system with a generalised coordinate $q$ and a bump potential at zero. The potential chosen is $V(q) = A\, \mathrm{sech}^2 (q)$ which is a reduced form of the Eckart potential barrier used to represent molecular reactions under stationary (gas-phase)  environments \cite{craven2017lagrangian}. This potential is symmetric about the origin where it has a global maximum. It follows that $V'(q) = -2 A\tanh (q) \mathrm{sech}^2 (q)$ and so the conservative system is given by the equation of motion.    
\begin{equation}
\label{eq: 1dof conservative}
    I \ddot{q}  = 2 A \tanh (q) \mathrm{sech}^2 (q),
\end{equation}
where $I >0$ corresponds to the mass of particle. We note that $H = \tfrac12 I\dot{q}^2 + V(q) $ is conserved under the flow, and the stable manifold to the equilibrium can be shown to be the graph given by $\dot{q} = -\sqrt{\frac{2A}{I}} \tanh q$ (see Figure \ref{fig:stableManifoldConservative}).

\subsubsection{Forced and damped model}
\label{subsec: 1dof forced damped model}
The model in Equation (\ref{eq: 1dof conservative}) is now perturbed by a linear damping term and a time-dependent forcing term, leading to the equation
\begin{equation}
    \label{eq: 1dof forced and damped}
    I \ddot{q}  = 2 A \tanh (q) \mathrm{sech}^2 (q) - k \dot{q} + F(t),
\end{equation}
where $k$ is a positive constant. Through the introduction of dimensionless  variables $s = \sqrt{\frac{2A}{I}}t$, $x(s) = q(t)$, and using $'$ to denote $\frac{\mathrm{d}}{\mathrm{d}s}$, Equation (\ref{eq: 1dof forced and damped}) can be recast as
\begin{equation}
    \label{eq: 1dof forced and damped rescaled}
    x'' = \tanh (x) \mathrm{sech}^2 (x) - \tilde k x' + \tilde F(s),
\end{equation}
where $\tilde k, \tilde F$ represent the re-scaling of the damping and forcing terms. Equation (\ref{eq: 1dof forced and damped rescaled}) will now be the object of the study of this subsection, as it represents the qualitative dynamics of Equation (\ref{eq: 1dof forced and damped}). From now on tildes will be omitted, and we will use $k,F$ to denote the dimensionless parameters.

To approach this system using the methods outlined in Section~\ref{sec: theory and methods}, we introduce $y = x'$ and the equations of motion become
\begin{equation}
    \label{eq: 1dof forced and damped rescaled, system}
    \begin{pmatrix}
        x'(s) \\
        y'(s)
    \end{pmatrix} =
    \begin{pmatrix}
        y(s) \\
        \tanh (x(s)) \mathrm{sech}^2 (x(s)) -  k y(s)
    \end{pmatrix}
    + \begin{pmatrix}
        0 \\
         F(s)
    \end{pmatrix}.
\end{equation}
In the case of $ F = 0$, Equation (\ref{eq: 1dof forced and damped rescaled, system}) has a unique hyperbolic steady state at $(0,0)$. Applying the terminology established in Section~\ref{sec: linear hyperbolic trajectory}, the linearised system corresponding to Equation (\ref{eq: 1dof forced and damped rescaled, system}) is given by 
\begin{equation}
    \label{eq: 1dof forced and damped rescaled, linearised system}
    \begin{pmatrix}
        u'(s) \\
        v'(s)
    \end{pmatrix} =
    \begin{pmatrix}
    0 & 1 \\
    1 & -k
    \end{pmatrix}
    \begin{pmatrix}
        u(s) \\
        v(s)
    \end{pmatrix}
    + \begin{pmatrix}
        0 \\
         F(s)
    \end{pmatrix},
\end{equation}
where $(u,v)$ denotes a small perturbation from $(0,0)$. This is a system of the form of Equation  \eqref{eq: linearised system} with
\begin{equation}\mathbf{A} =    \begin{pmatrix}
    0 & 1 \\
    1 & -k
    \end{pmatrix}.  \end{equation}
The origin is therefore a saddle point with eigenvalues $\lambda_\pm = \frac{- k \pm \sqrt{4 +  k^2}}{2}$ and eigenvectors $\mathbf{v}_\pm = (1,\lambda_\pm)^T$. Letting $\mathbf{P} = (\mathbf{v}_-,\mathbf{v}_+)$ be the change of basis matrix to the eigenvectors, it can be shown that the system defined by Equation  \eqref{eq: 1dof forced and damped rescaled, linearised system} is equivalent to the system:
\begin{equation}
    \label{eq: 1dof forced and damped rescaled, linearised and transformed system}
    \begin{pmatrix}
        y_-'(s) \\
        y_+'(s)
    \end{pmatrix} =
   \begin{pmatrix}
    -1 - \lambda_-^2 & 0\\
   0& 1 + \lambda_+^2 
    \end{pmatrix}
    \begin{pmatrix}
        y_-(s) \\
        y_+(s)
    \end{pmatrix}
    + \begin{pmatrix}
        - F(s) \\
         F(s)
    \end{pmatrix},
\end{equation}
where $(u,v)^T = y_- \mathbf{v}_- + y_+  \mathbf{v}_+$. Equation  \eqref{eq: 1dof forced and damped rescaled, linearised and transformed system} has the same form
as the two-dimensional system in Section \ref{sec: linear hyperbolic trajectory}, and so can be solved using the general solution presented there and then transformed to the original coordinate system.
\subsubsection{Results}
\label{sec: 1dof results}
\begin{figure}[hbtp]
    \begin{subfigure}[b]{ 0.48\textwidth}
    \centering
    \includegraphics[width = \textwidth]{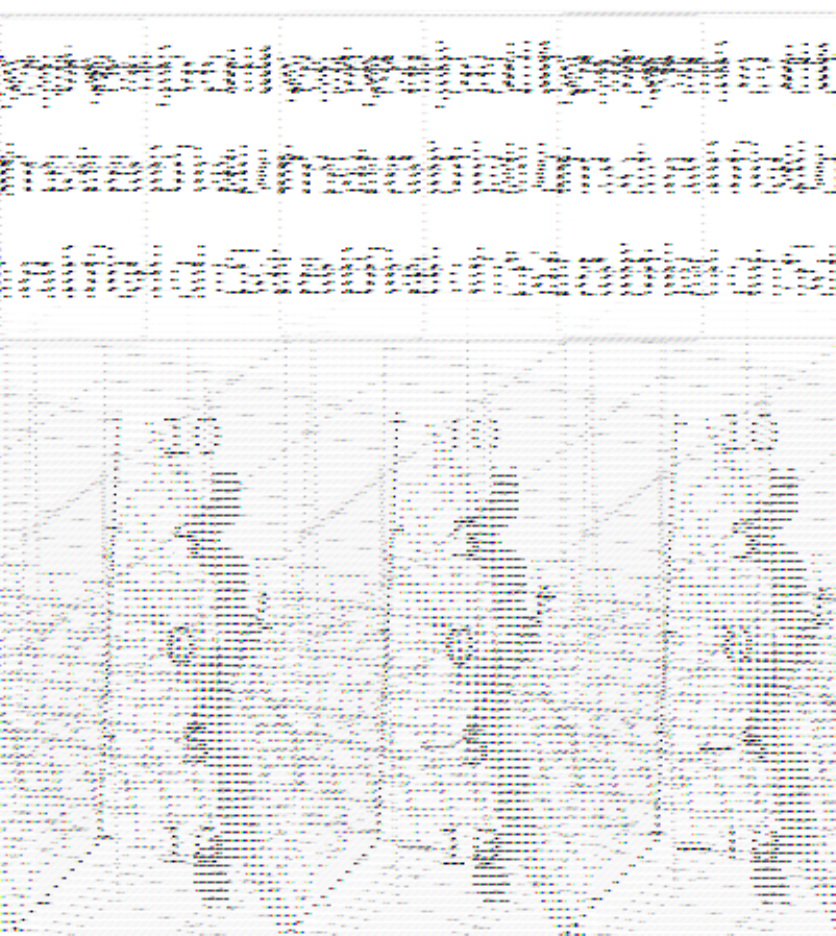}
    \caption{$k=1$}
    \end{subfigure}
    \begin{subfigure}[b]{ 0.48\textwidth}
    \centering
    \includegraphics[width = \textwidth]{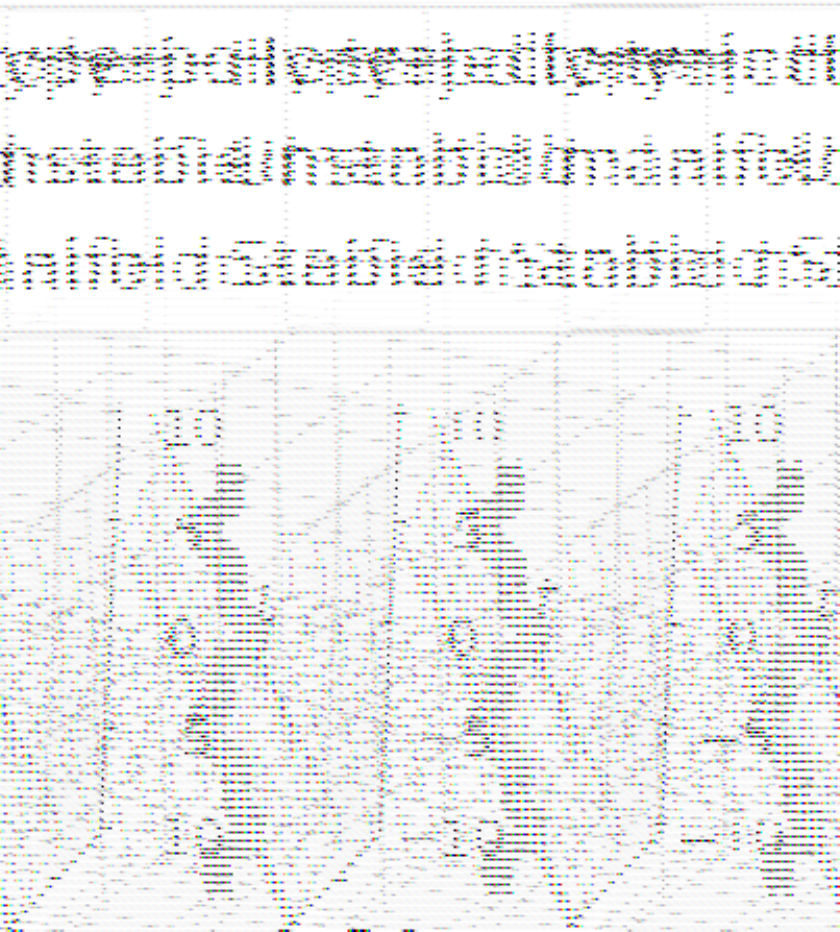}
    \caption{$k=3$}
    \end{subfigure}
    \caption{The stable and unstable manifolds to the hyperbolic trajectory of $x'' = \tanh (x) \mathrm{sech}^2 (x) - k x' + 0.1 \cos(0.7s) -0.15 \cos(2s)$ for $k = 1,3$.}
    \label{fig: 1dof stable unstable original coordds}
\end{figure}

The methods outlined in Sections \ref{sec: algorithm hyp traj} and \ref{subsec: algorithm stable manifold} were applied to the system in Equation  \eqref{eq: 1dof forced and damped rescaled, system}, using the linearised system given by \eqref{eq: 1dof forced and damped rescaled, linearised and transformed system} for the initial guess of the trajectories. For numerical simulations we choose $k = 1$ or $3$, and $ F(s) = 0.1 \cos(0.7s) -0.15 \cos(2s)$. The time interval was chosen to be $[-10,10]$ and the discretisation parameter $N$ was chosen to be 401.
The convergence parameters were set to be $\varepsilon_F = 10^{-6}$ and $\varepsilon_C=10^{-7}$. For both values of $k$ the hyperbolic trajectory was obtained in 7 iterations of the Newton-Raphson method.
In Figure \ref{fig: 1dof stable unstable original coordds} a hyperbolic trajectory for this forcing (and for two values of damping) is plotted together with the stable and unstable manifolds in the extended state space of $\mathbb{R}^2 \times \mathbb{R}$. 

\begin{figure}[hbtp]
    \begin{subfigure}[b]{0.48\textwidth}
    \centering
    \includegraphics[width = \textwidth]{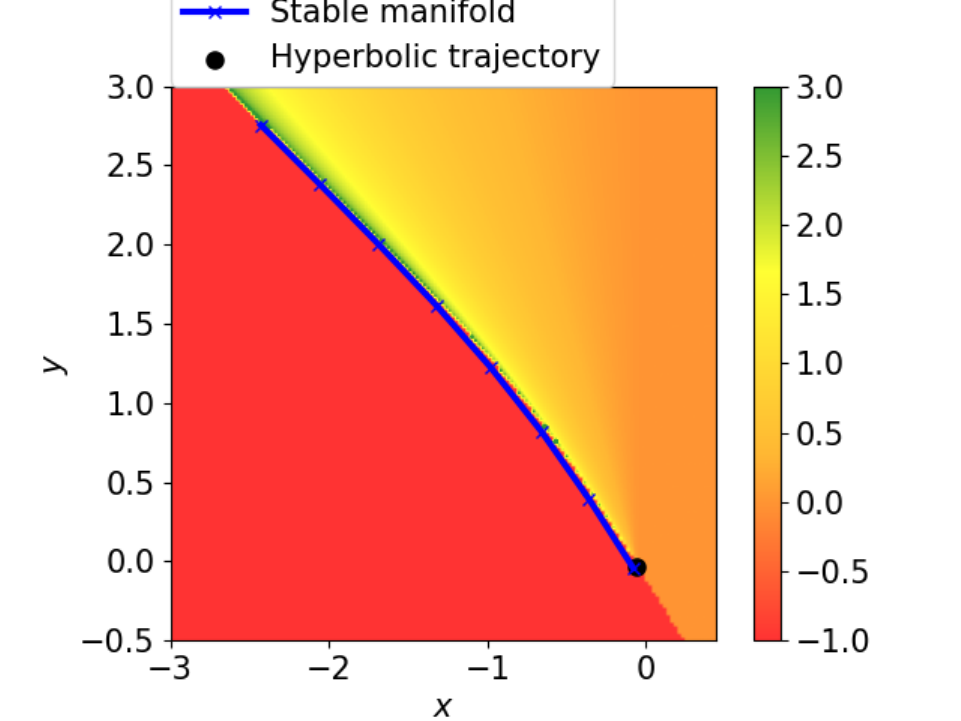}
    \caption{$k=1$}
    \end{subfigure} \quad
    \begin{subfigure}[b]{0.48\textwidth}
    \centering
     \includegraphics[width = \textwidth]{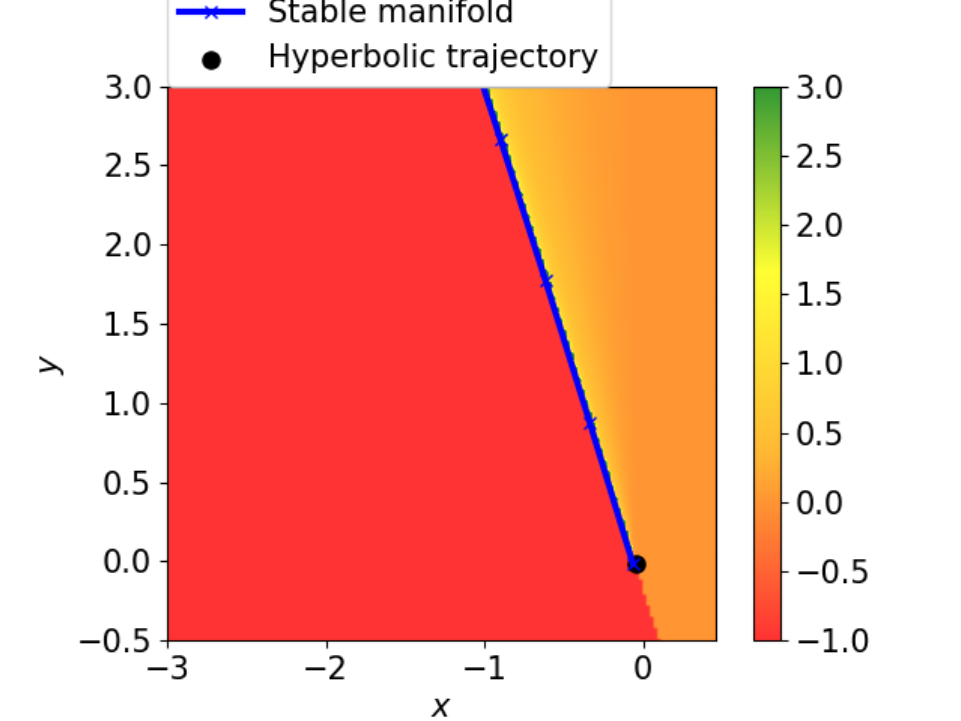}
     \caption{$k=3$}
    \end{subfigure}
    \caption{A branch of the stable manifold is plotted in the plane $t=0$ and this is compared with a heat map of initial conditions coloured according to the time it takes to transition over the hyperbolic trajectory. The points which do not transition are assigned the value of $-1$ and appear as red. This is done for the systems $x'' = \tanh (x) \mathrm{sech}^2 (x) - k x' + 0.1 \cos(0.7s) -0.15 \cos(2s)$ with $k = 1,3$.}
    \label{fig: 1dof stable unstable timeslice}
\end{figure}

A trajectory is deemed to cross the hyperbolic trajectory when it crosses the dividing manifold given by  $x=x_{hyp}(t)$.  This manifold consists of two parts: $y > y_{hyp}(t)$, which trajectories cross with $\dot{x}>0$; and $y < y_{hyp}(t)$, which trajectories cross with $\dot{x}<0$; this allows us to compare the direction of the crossing of a trajectory. In Figure \ref{fig: 1dof stable unstable timeslice} the stable manifold is presented along with a heat map of the transition time, which was measured as the first time a solution crosses the line $y = x_{hyp}(s)$. The stable manifold consists of states which converge to the hyperbolic trajectory but will not cross over it, and this is reflected in the states close to the stable manifold in the transition region having a greater transition time.
\subsubsection{Accuracy}
The standard method to compute invariant manifolds of hyperbolic trajectories is to begin with a small initial approximation of the manifold, and evolve this approximation in time, inserting additional points when required \cite{mancho2003computation}. The case of hyperbolic trajectories in a two-dimensional system subject to aperiodic forcing was studied by \cite{mancho2003computation} and we shall use this as a benchmark to test our method.

\begin{figure}[hbtp]
    \begin{subfigure}[b]{ 0.45\textwidth}
    \centering
    \includegraphics[width = \textwidth]{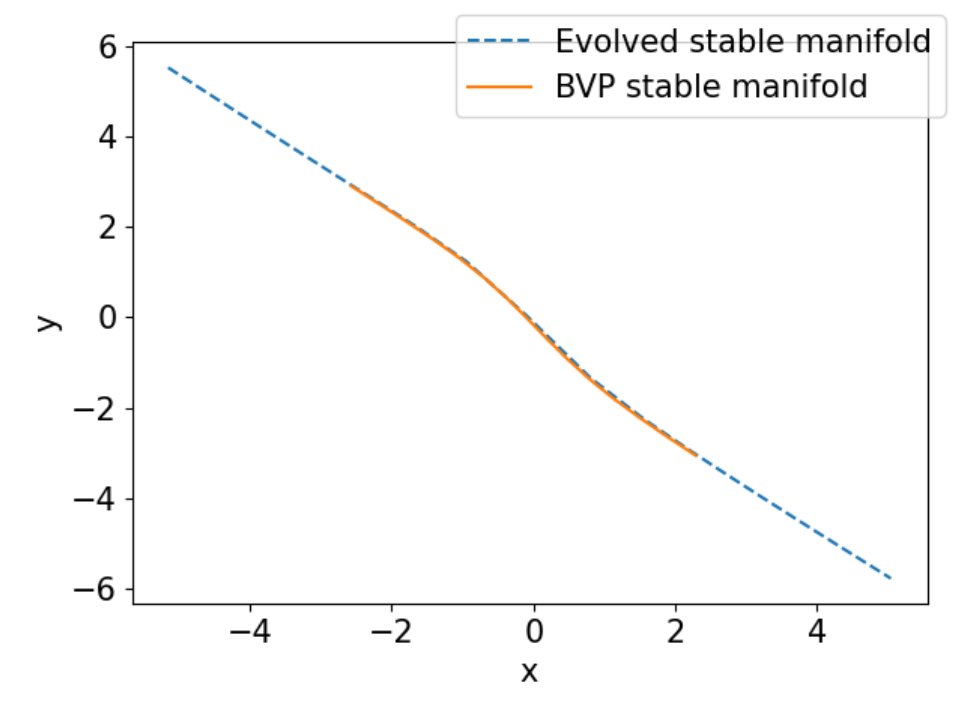}
    \caption{Stable manifold for $k=1$}
    \end{subfigure}
    \begin{subfigure}[b]{ 0.45\textwidth}
    \centering
    \includegraphics[width = \textwidth]{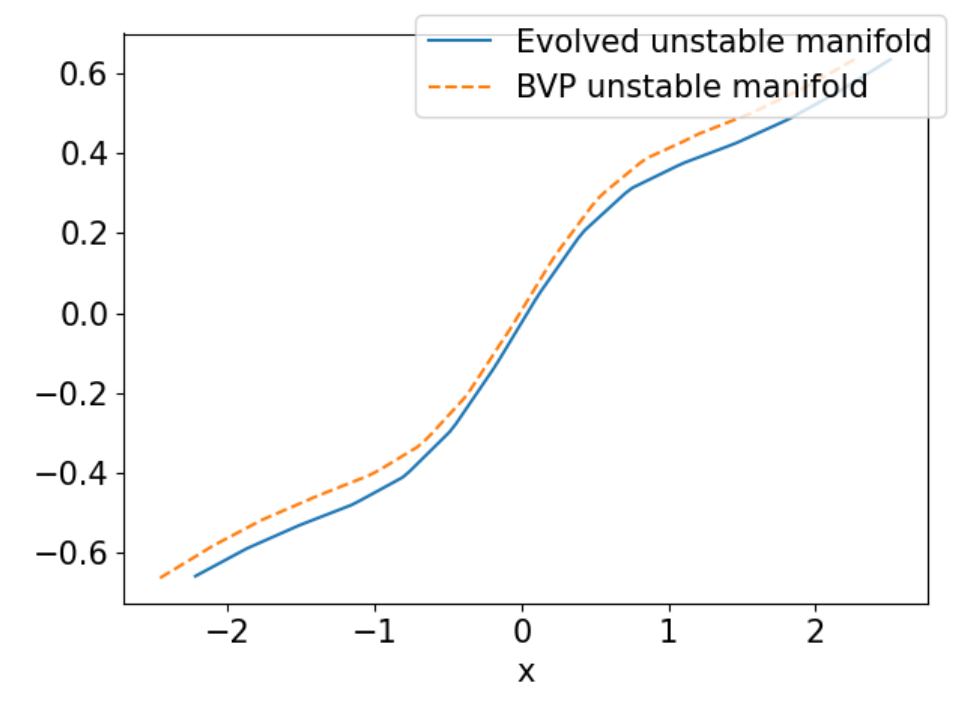}
    \caption{Unstable manifold for $k=1$}
    \end{subfigure}
    \begin{subfigure}[b]{ 0.45\textwidth}
    \centering
    \includegraphics[width = \textwidth]{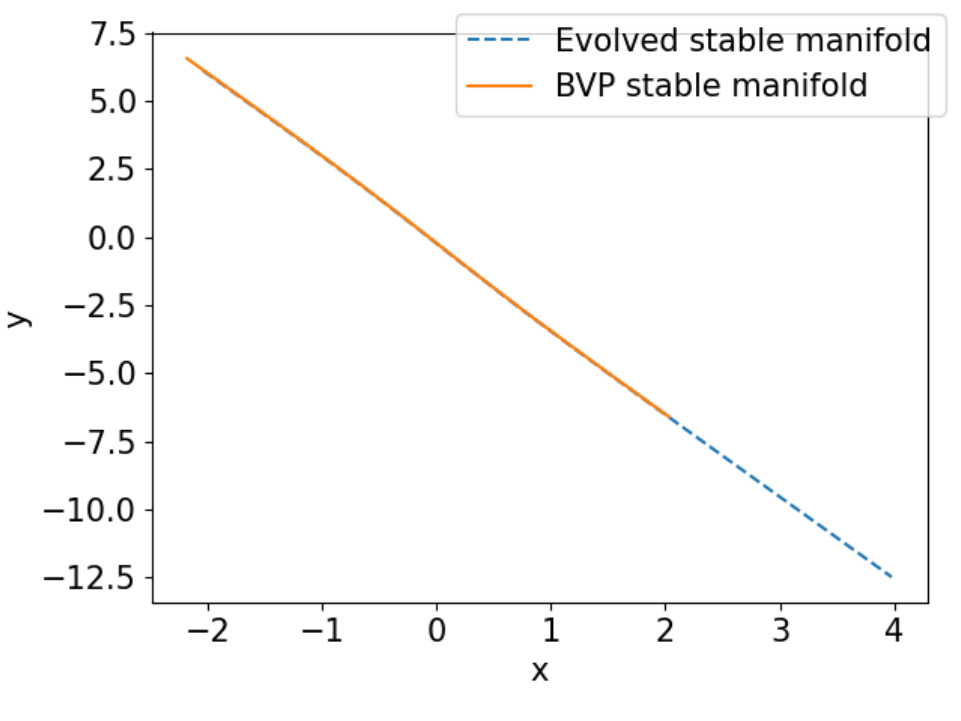}
    \caption{Stable manifold for $k=3$}
    \end{subfigure}
    \begin{subfigure}[b]{ 0.45\textwidth}
    \centering
    \includegraphics[width = \textwidth]{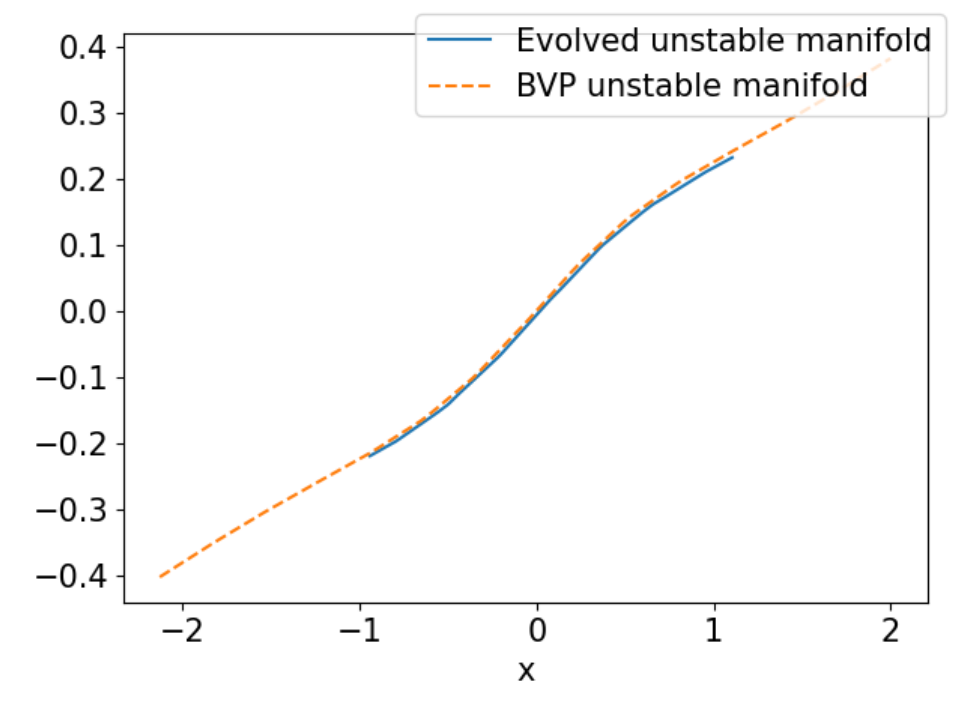}
    \caption{Unstable manifold for $k=3$}
    \end{subfigure}
\caption{Comparison of manifolds obtained using methods from Section \ref{sec: theory and methods} and from evolving an initial segment of the manifolds. Using the method from Section \ref{sec: theory and methods}, 100 points were initially sampled along each of the manifolds at times 
$t = -5.0, 1.75$ for the unstable and stable manifolds with $k=3$ and $t = -6.0, 4.25$ for the unstable and stable manifolds with $k=1$. These were evolved backwards/forwards to $t = 0$, using the Hobson criteria \cite{mancho2003computation} with $(\alpha,\mathrm{d}\alpha,\delta) = (0.3,0.0001,0.000001)$ and $\mathrm{d}T = 0.005$. Points were inserted using cubic splines to interpolate the manifolds.}
    \label{fig: bvp compared with Mancho}
\end{figure}

In Figure \ref{fig: bvp compared with Mancho} we can see that our method does accurately represent the invariant manifolds at $t = 0$. We first computed using the methods from Section \ref{sec: theory and methods} the stable manifolds for the system $t = 0$, with $k = 1, F(t) = 0.1\cos(0.7t) - 0.15\cos(2t)$. One hundred points were sampled along each manifold, at times 
t = -5.0,1.75 for the unstable and stable manifolds with $k=3$ and t = -6.0,4.25 for the unstable and stable manifolds with $k=1$. These branches of the stable and unstable manifold were evolved forwards and backwards respectively to $t=0$ with a step size of $dT = 0.005$. The initial time values for the stable and unstable manifolds were chosen experimentally to ensure that when evolving these sections forwards/backwards in time to $t=0$, the lengths of the sections are similar to the length of the manifolds computed directly from the BVP method at $t=0$. To determine when it was necessary to add points into the discretisation we used the Hobson criteria established in \cite{mancho2003computation} and cubic splines were used to insert the points.

\subsection{Two DoF model}
\label{subsec: two dof model}
In this section, the methods outlined in Section \ref{sec: theory and methods} are applied to  a qualitative roll-heave model  from \cite{thompson1996suppression}, to model ship motions. In contrast to the analysis done in \cite{thompson1996suppression}, which focused on a roll-heave system with periodic forcing, our method extends to arbitrary bounded forcing functions. In \cite{thompson1996suppression} the attractor of the system could be studied through the behaviour of the Poincar\'e mapping, which is not possible for non-periodic forcing. This makes the extended state-space 5-dimensional, which is harder to visualise, but it captures the dynamics of escape better.

In the absence of damping or external forces acting on the ship, the system is taken to evolve according to
\begin{align}
m\ddot{z} +h\left( z - \tfrac12 \gamma \phi^2 \right) & = 0 \\
I \ddot{\phi} +c \phi\left[1-\left(\phi / \phi_{v}\right)\right]\left[1+\left(\phi / \phi_{v}\right)\right]  - h \gamma \phi\left( z - \tfrac12 \gamma \phi^2 \right) &= 0.
\label{eq: thompson  equations of motion}
\end{align}
To this equation we will add forcing and linearised damping. We then study a dimensionless form of the equation, by rescaling the $z,\phi,t$ to dimensionless quantities $x,y,s$, so that Equation \eqref{eq: thompson  equations of motion} is equivalent to 
\begin{equation}
(x',y',v_x',v_y')^T =\mathbf{f}(x,y,v_x,v_y) \label{eq: 2dof rescaled vector}
\end{equation}
where 
\begin{equation}
\label{eq: 2dof rescaled f}
\mathbf{f}(x,y,v_x,v_y) = (
    v_x ,    v_y,     -  h (x-y^2) - k_x v_x ,    -y + xy- k_yv_y)^T.
\end{equation}
The details of this rescaling are presented in Appendix \ref{appendix: derivation}, and the analysis of the equilibria of this system is in Appendix \ref{appedix: analysis}. In Appendix \ref{appendix:autoHam} we present for comparison the analysis of flux over the saddles in the autonomous Hamiltonian case.

As in the previous example, the framework established in Section \ref{sec: theory and methods} can be applied to the system in Equation~\eqref{eq: 2dof rescaled vector}. For this system there  are two saddle points corresponding to $y = \pm 1, x= 1, v_x = v_y = 0$  and hence we will obtain two hyperbolic trajectories and two stable manifolds. We shall refer to the hyperbolic trajectory continued from the saddle point $(1,1,0,0)$ as the \emph{positive hyperbolic trajectory}  and the hyperbolic trajectory continued from $(1,-1,0,0)$ as the  \emph{negative hyperbolic trajectory}.

\subsubsection{Numerical experiments}
\label{sec: numerical experiments}
\begin{figure}
    \centering
    \begin{subfigure}[b]{.45\textwidth}
        \includegraphics[width=\textwidth]{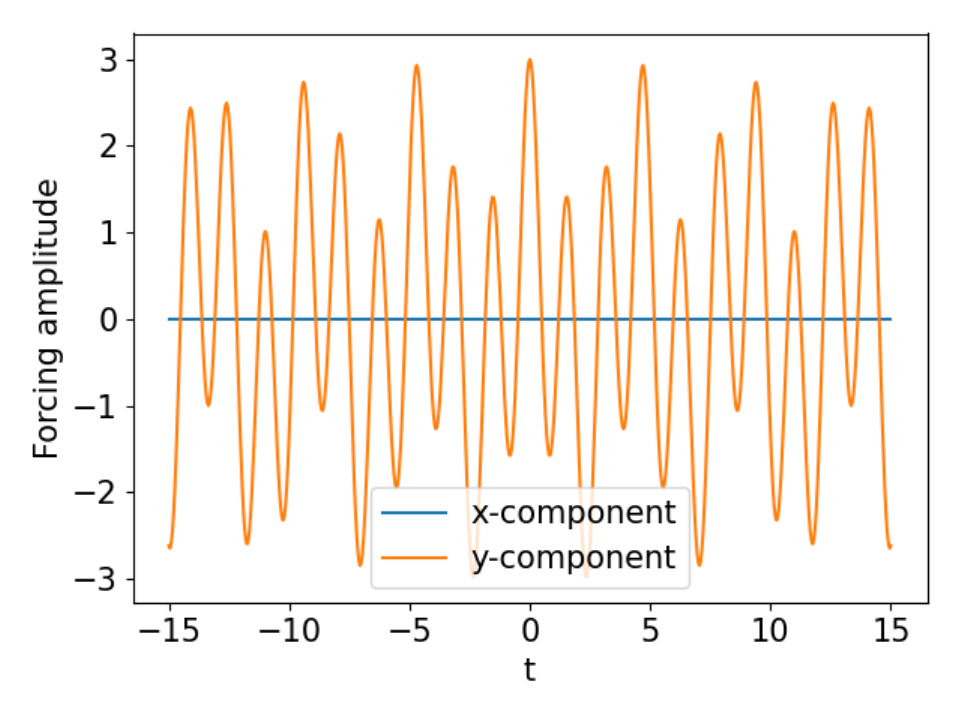}
        \caption{Forcing components for quasi-periodic forcing, $F(t)=0,M(t) = \cos(\sqrt{2}t) + 2\cos(4t)$}
        \label{fig:quasi-periodic forcing}
    \end{subfigure} \quad
        \begin{subfigure}[b]{.45\textwidth}
        \includegraphics[width=\textwidth]{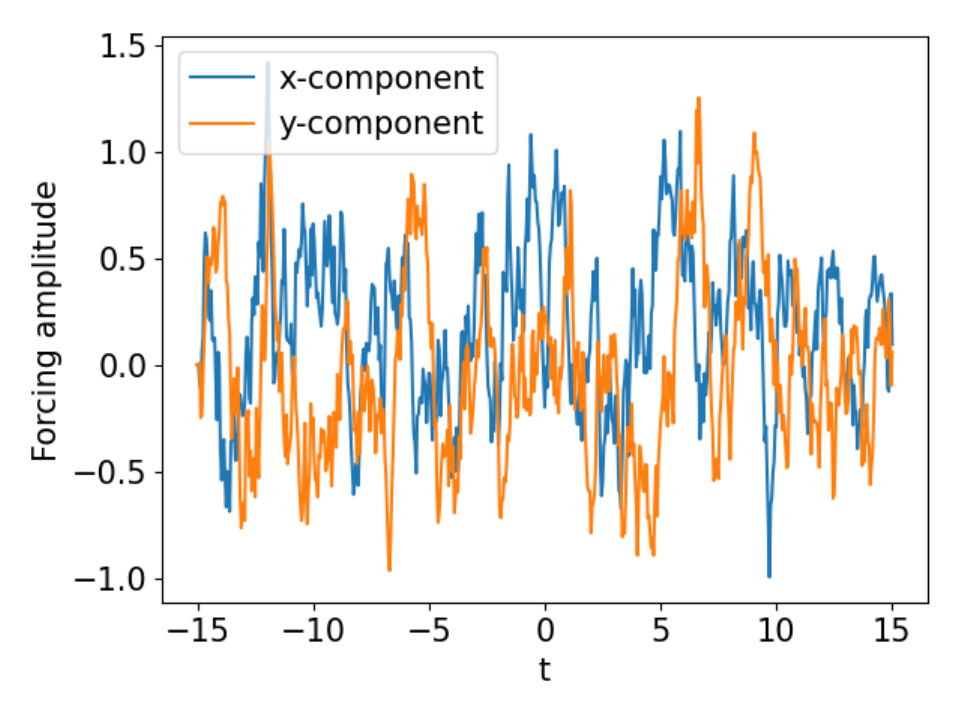}
        \caption{Forcing components for filtered white noise forcing, $\eta(t) = (F(t),M(t))^T$}
        \label{fig:stochastic forcing}
    \end{subfigure}
    \caption{Forcing components for two degree of freedom system.}
    \label{fig:enter-label}
\end{figure}

The methods described in this paper were applied to the system  \eqref{eq: 2dof rescaled system} with  model parameters $h=k_x=k_y=1$ and two cases of external forcing. These forcing functions are a quasi-periodic forcing function in the roll-direction $F(t)=0,M(t) = \cos(\sqrt{2}t) + 2\cos(4t)$ (Figure \ref{fig:quasi-periodic forcing}) and a filtered white noise force $\eta(t) = (F(t),M(t))^T$  (Figure \ref{fig:stochastic forcing}), a sample from the stochastic process
\begin{equation}
\label{eq:stochastic forcing}
    d\eta = \begin{pmatrix}-\lambda & -\omega \\ \omega & -\lambda\end{pmatrix} \eta \, dt + B\, dW.
\end{equation}
Here, $dW$ represents the increments of a standard 2D Wiener process. 
We choose $\omega =\tfrac12 \sqrt{\alpha - k^2}$, to match  the angular frequency of oscillation about the saddle as a likely worst case, and $\lambda = \lambda_3$, the strength of contraction along the strong stable manifold (for no particular reason; $\lambda$ represents the autocorrelation decay rate of the filtered white noise). We choose $B = \tfrac12 \begin{pmatrix}
    1 & 1 \\ -1 & 1
\end{pmatrix}$. A sample solution of equation \eqref{eq:stochastic forcing} was found by computing the Ito integral using the python library sdeint \cite{sdeint}, but it can be done explicitly, e.g.~\cite{mackay2021inference}.

Now that the system has been established our goals are to: compute the hyperbolic trajectories, compute the stable manifolds, compute a dividing manifold, and use these to classify initial states and their time to capsize. The details and the results of these experiments are in the following section.

\subsubsection{Results}
\label{subsec: ship results}
The time-interval was chosen to be the symmetric interval $[-15,15]$ ($T=15$) and the number of sample times $N$ chosen to be 601 (an odd number so that $t=0$ is a sample time). 

This system has two hyperbolic trajectories of interest, corresponding to the saddle points $(1,\pm1,0,0)$. These were first computed according to the method described in Section \ref{sec: algorithm hyp traj}. For the initial choice of the hyperbolic trajectory we chose the constant function ($t\mapsto (1,\pm1,0,0)$) and the convergence parameters were set to be $\varepsilon_c=10^{-5},\varepsilon_f=10^{-6}$.
Table \ref{tab: iterations to hyp traj} demonstrates the convergence of the method to the hyperbolic trajectory.
\begin{table}[ht!]
    \centering
    \begin{tabularx}{\textwidth}{|c|X|X|}
    \hline
    \textbf{Forcing} & \textbf{Negative trajectory\newline iterations}  & \textbf{Positive trajectory \newline itertations}  \\
    \hline
    Quasi-periodic    &8 &8 \\
    Filtered white noise & 12&12\\
    \hline
    \end{tabularx}
    \vspace*{5mm}
    \caption{Iterations to find hyperbolic trajectories for the system given by \eqref{eq: 2dof rescaled system}, using the method outlined in Section \ref{sec: algorithm hyp traj}.}
    \label{tab: iterations to hyp traj}
\end{table}
The hyperbolic trajectories corresponding to the saddle points $(1,\pm1,0,0)$ are plotted in $(x,y,t)$-space in Figure~\ref{fig: hyp trajs in x,y space}.
\begin{figure}[hbtp]
    \centering
    \begin{subfigure}[b]{.49\textwidth}
        \includegraphics[width=\textwidth]{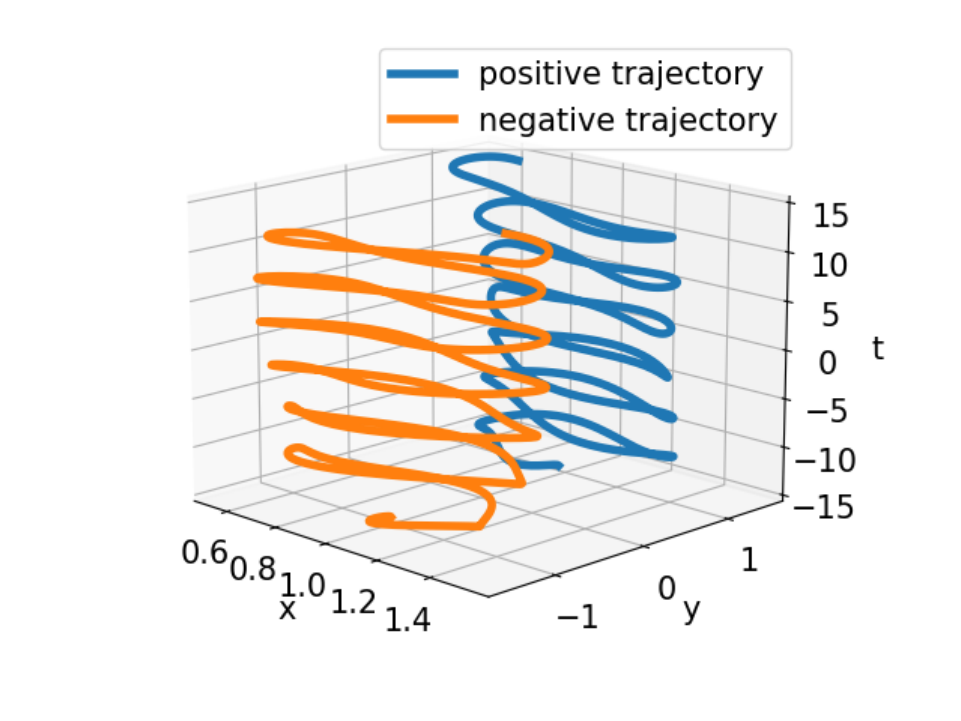}
        \caption{Quasi-periodic forcing}
    \end{subfigure}
    \begin{subfigure}[b]{.49\textwidth}
        \includegraphics[width=\textwidth]{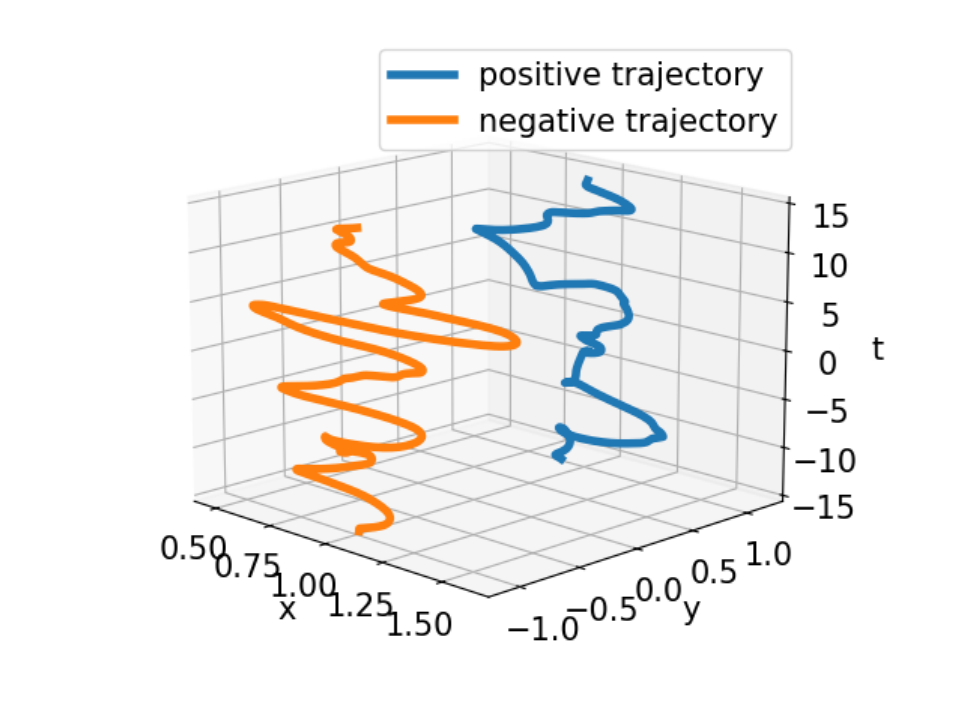}
        \caption{Filtered white noise forcing}
    \end{subfigure}
    \caption{Plots of the hyperbolic trajectories in $(x,y,t)$-space for quasi-periodic and filtered white noise forcing.}
    \label{fig: hyp trajs in x,y space}
\end{figure}

\begin{figure}[hbtp]
    \centering
    \begin{subfigure}[b]{.99\textwidth}
    \includegraphics[width=\textwidth]{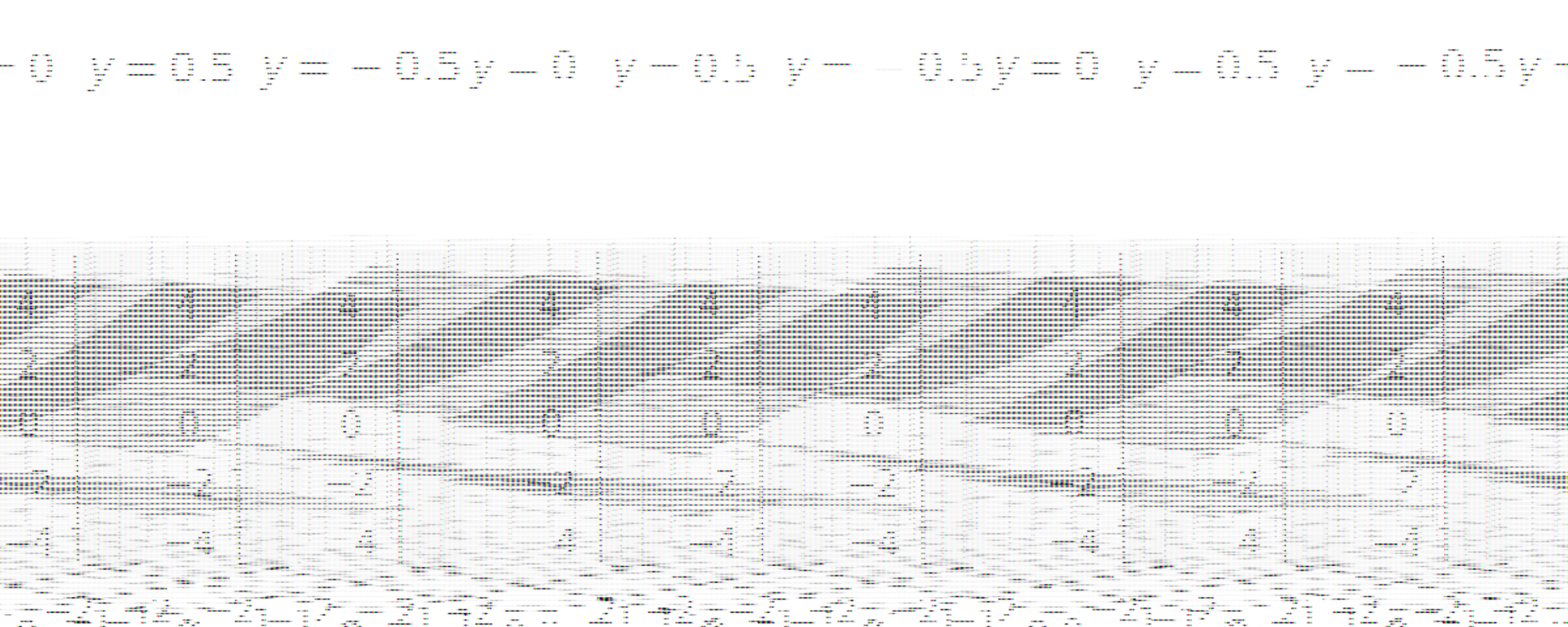}
    \caption{Quasi-periodic forcing}
    \end{subfigure}
    \begin{subfigure}[b]{.99\textwidth}
    \includegraphics[width=\textwidth]{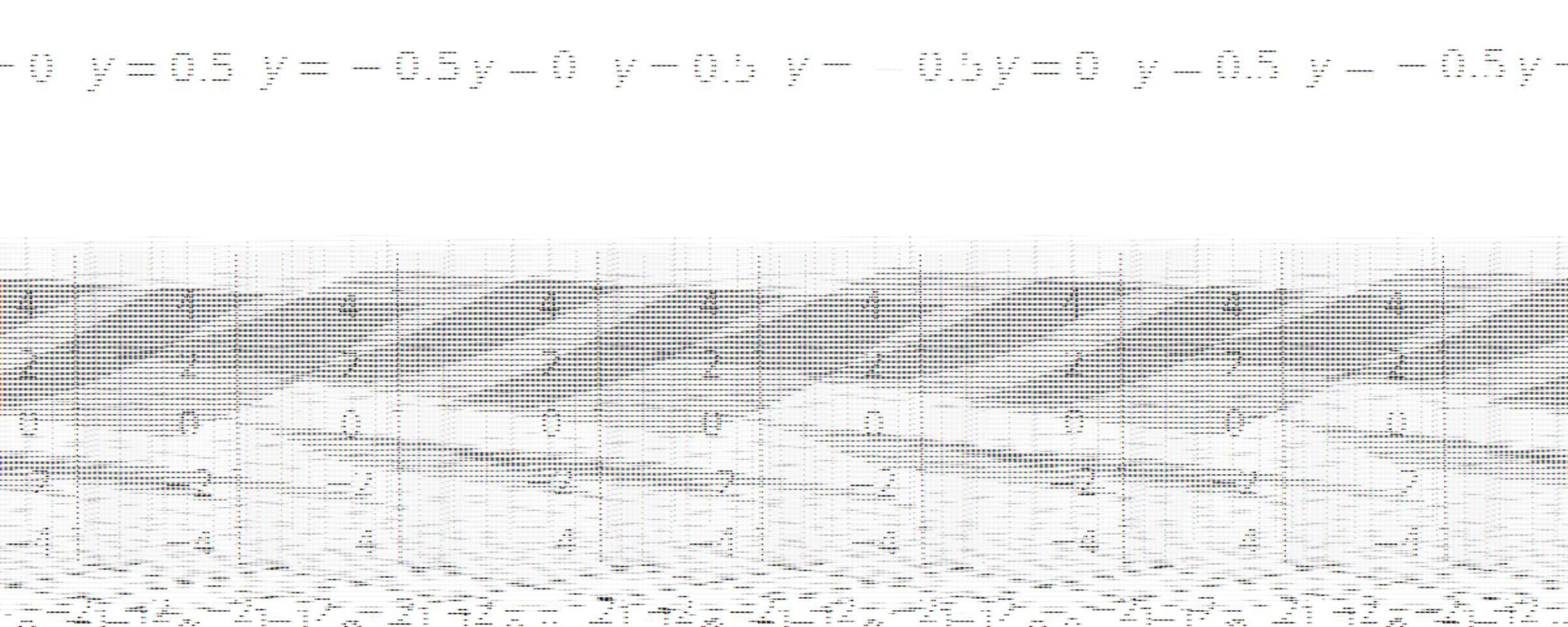}
    \caption{Filtered white noise forcing}
    \end{subfigure}
    \caption{Plots of the stable manifolds as a graph of the $v_y$ component over the $(x,y,v_x)$ components. For each plot, a value of $y$ is fixed and $(x,v_x)$ are allowed to vary.}
    \label{fig:vy graphs}
\end{figure}

Similarly the stable and centre manifolds for these hyperbolic trajectories were computed according to the method described in Section \ref{sec: algorithm center manifold} and \ref{subsec: algorithm stable manifold}. First we computed the stable manifold on the time slice $t=0$. This was done by applying the method in Section \ref{subsec: algorithm stable manifold} for $j=301$. The boundary value problem parameter $\mathbf{q}$ was sampled uniformly from a hypercube lattice in $[-1.5,1.5]^3$. \\
Upon observation it was found that these stable manifolds are given locally by graphs of the $v_y$ variable as a function $\tilde{v}_y$ of the $(t,x,y,v_x)$ variables.   Plots of these are given in Figure \ref{fig:vy graphs}. 

Using a multiquadric radial basis function interpolation (see \cite{buhmann2000radial} for background), the function $\tilde v_y(t,x,y,v_x)$ was approximated, and for a given  initial condition $(\mathbf{x},t) = (x,y,v_x,v_y,t)$ it was possible to classify it as capsize or non-capsize by considering the sign of $v_y - \tilde v_y(t,x,y,v_x)$. This must be done for the stable manifold corresponding to each of the centre manifolds for each of the saddle points $(1,\pm 1,0 ,0)$. 

We tested this classification procedure by direct integration of the differential equation for 10000 points randomly sampled from the hypercube $[-1,1]\times[-1,1]\times[-5,5]\times[-5,5]$. For the integration method, we considered a point to lead to capsize if the $y$ component escapes the region $\{y^2 <10\}$ while $t\le t_{max} = 11.5$. These thresholds were chosen from observing simulations and noticing that they are sufficiently large in the context of the system to ensure that capsize will occur.  Accuracy was compared by the ratio of agreement between the two methods. We found that for the quasi-periodic forcing the accuracy was 96.9\%, for the filtered white-noise it was 93.6\% .   The sensitivity (probability a capsize classification is correct) was 92.7\% for quasi-periodic forcing and 96.2\% for filtered white-noise forcing. The specificity (probability a safe classification is correct) was 99.8\% for quasi-periodic forcing and 92.0\% for filtered white-noise forcing.

Visualising this is difficult as even in this reduced time-slice $t=0$ the structure is inherently 4-dimensional. To present visually the distinction between capsize and non-capsize region we sampled points uniformly at random on the planes $t=0,v_x=0,y= -0.5,0,0.5$, that is we have fixed $t,v_x,y$ and vary $(x,v_y)$. The results are shown in Figures \ref{fig: scatter plot in slices, quasi-periodic} and \ref{fig: scatter plot in slices, stochastic}.

We note here that unlike in \cite{thompson1996suppression} we do not see the folding of the stable manifolds onto themselves. This is because our method is only valid for the region of the stable manifold which is given locally by graph over an affine subspace of the extended state space. We expect that if these sections of manifolds were integrated backwards in time we would observe such folding of the stable manifold. In practice however, the computational complexity of this task is significant, as numerical integration of the manifolds must also be combined with insertion of additional points.
\begin{figure}[hbtp]
    \centering
    \includegraphics[width=\textwidth]{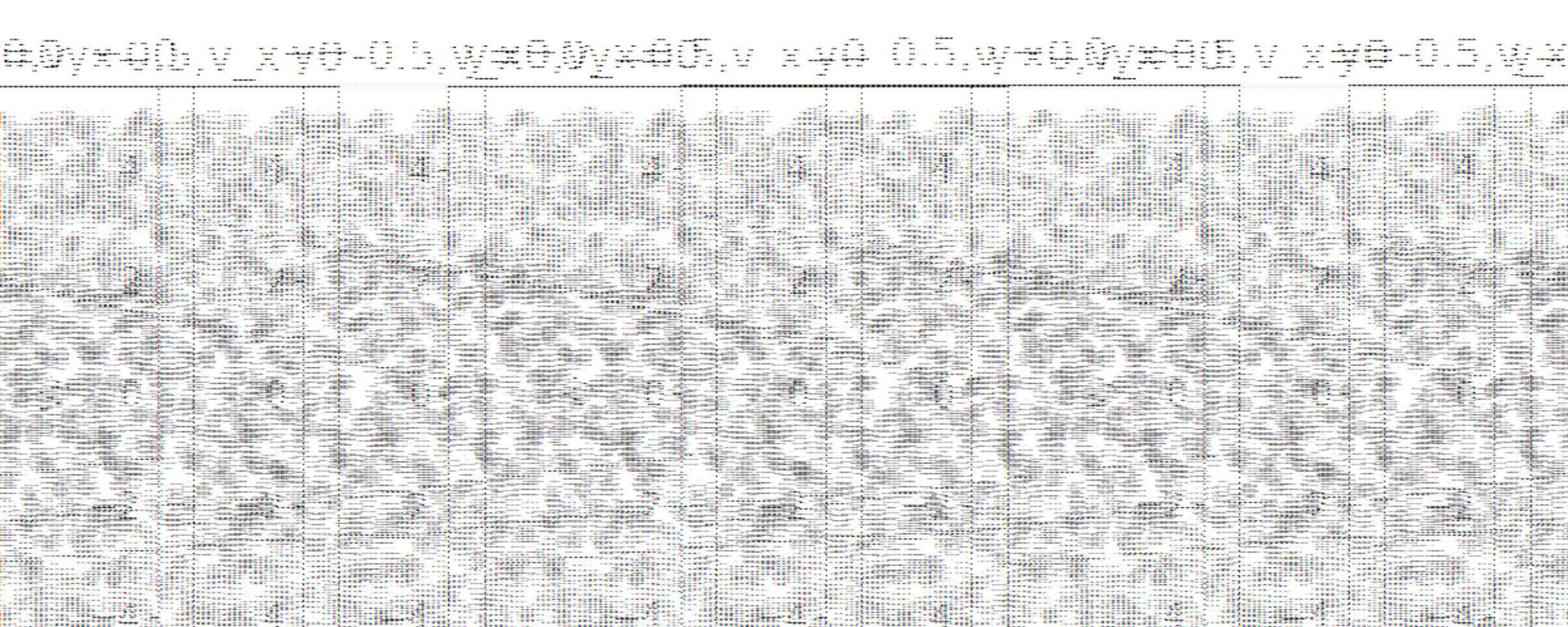}
    \caption{Scatter plot of capsize points for the system with quasi-periodic forcing in the $t=0,v_x=0,y=$const planes, determined by integration of the system and how they lie between the stable manifolds.}
    \label{fig: scatter plot in slices, quasi-periodic}
\end{figure}
\begin{figure}[hbtp]
    \centering
    \includegraphics[width=\textwidth]{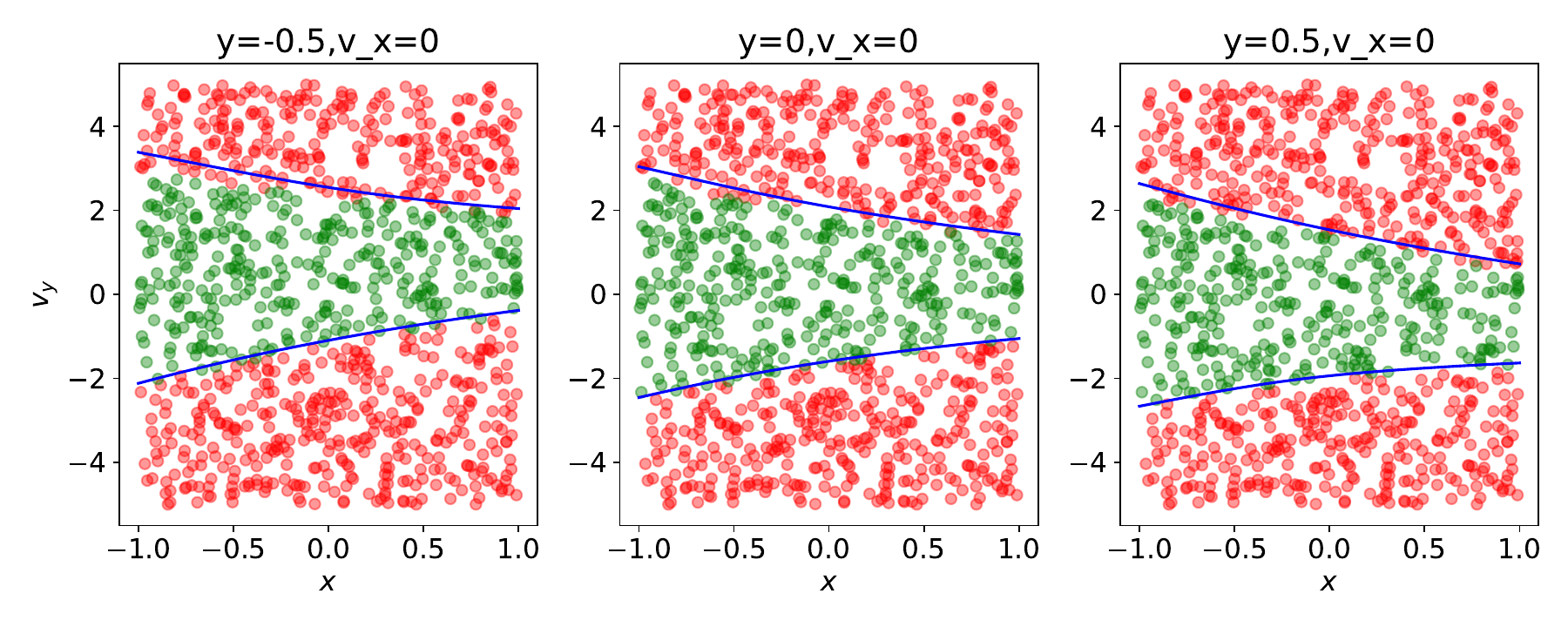}
    \caption{Scatter plot of capsize points for the system with filtered white noise forcing in the $t=0,v_x=0,y=$ const planes, determined by integration of the system and how they lie between the stable manifolds.}
    \label{fig: scatter plot in slices, stochastic}
\end{figure}

Once the capsize points were established we used a dividing manifold to define the time to capsize. We consider crossing the dividing manifold in the appropriate direction as implying capsize.  There is considerable arbitrariness to the dividing manifold, but they differ only in the time along a trajectory at which capsize is declared to take place. 
For our choice of dividing manifold corresponding to capsizing over $(1,\pm1,0,0)$ we took:
\begin{equation}
    \label{eq: dividing manifold}
    D_{M,\pm} = \{(x,y,v_x,v_y,t) \in \mathbb{R}^5 : (x,y,v_x,v_y) - \mathbf{x}_{hyp,\pm}(t) \in \Pi_{t,\pm}\},
\end{equation}
where $\mathbf{x}_{hyp,\pm}$ is the hyperbolic trajectory corresponding to the side of interest and
\begin{equation}\Pi_{t,\pm} = \mathrm{Span}(\mathbf{D}\iota_{S,\pm}(e_1),\mathbf{D}\iota_{S,\pm}(e_2),\tfrac12 \mathbf{D}\iota_{S,\pm}(e_3) + \tfrac12 \mathbf{v}_4 ), \end{equation}
where $\mathbf{v}_4$ corresponds to the unstable eigenvector direction. Equation  \eqref{eq: dividing manifold} can be interpreted as a smoothly varying family of affine subspace approximations to a dividing manifold, and so at each time it is straightforward to check which side of the dividing manifold a state is on, and hence it is simple to check when a trajectory does cross the dividing manifold, by expressing $D_{M,\pm}$ as $\{(\mathbf{x},t): L_t(\mathbf{x}) = 0 \}$, where for each $t$, $L_t$ is an affine linear function. 

The dividing manifold divides the state space into a \textit{capsized} and a \textit{not yet capsized} region, and the time to capsize of a trajectory is defined to be the time at which a trajectory crosses from one region into the other. By considering the sign of $L_t(\mathbf{x})$ for points known to lie in the capsize or non-capsize regions, these sides can be consistently categorised. 

To compute this requires a knowledge of the partial derivatives of the stable manifolds immersions for each value of $t$. These partial derivatives only need to be taken at the origin, $\mathbf{q}=\mathbf{0}$. Therefore for the values of $t$ in the discretised interval greater than $0$, the boundary value problem was solved for 125 points uniformly sampled from a lattice in $[-.5,.5]^3$. At each time sample the immersion $\iota_S$ was computed using a multiquadric  radial basis function interpolation of the points from these boundary value problems. 

We can visualise these for the points sampled in the $x,v_x$ planes. In Figures~\ref{fig: time plot in slices new, quasi-periodic} and \ref{fig: time plot in slices new, stochastic}, we colour each capsize point according to the capsize time. We can see how this decreases further away from the stable manifolds. 
\begin{figure}[hbtp]
    \centering
    \includegraphics[width=\textwidth]{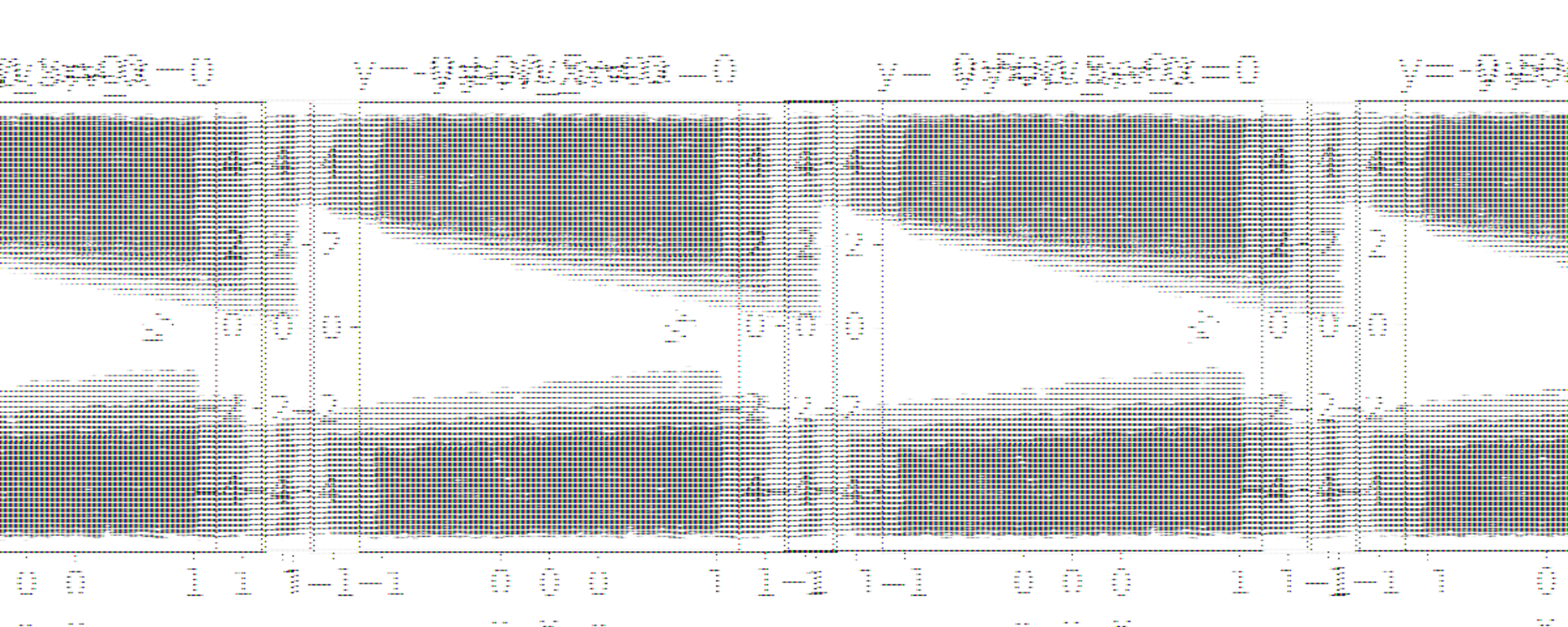}
    \caption{Scatter plot of capsize points coloured according to capsize time for the system with quasi-periodic forcing in the $t=0,v_x=0,y=$ const planes, determined by integration of the system and how they lie between the stable manifolds.}
    \label{fig: time plot in slices new, quasi-periodic}
\end{figure}
\begin{figure}[hbtp]
    \centering
    \includegraphics[width=\textwidth]{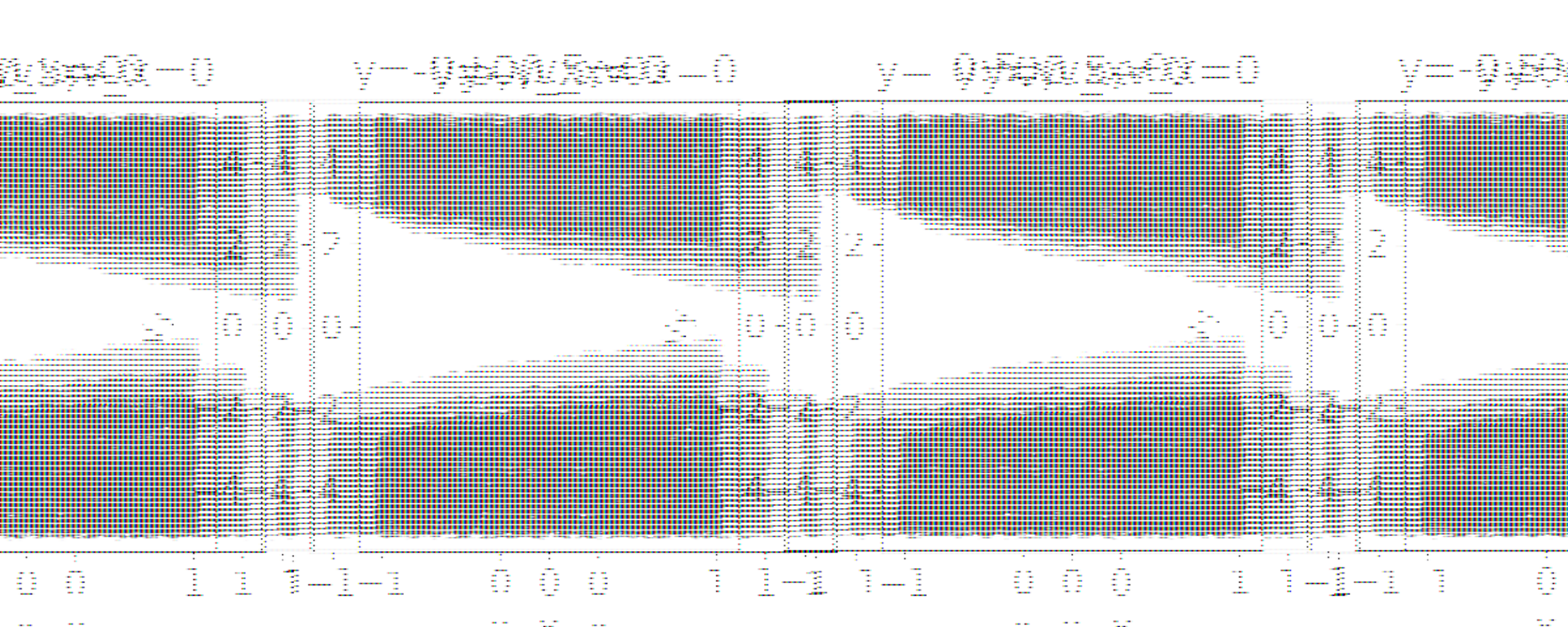}
    \caption{Scatter plot of capsize points coloured according to capsize time for the system with filtered white noise forcing in the $t=0,v_x=0,y=$ const planes, determined by integration of the system and how they lie between the stable manifolds.}
    \label{fig: time plot in slices new, stochastic}
\end{figure}

\subsubsection{Integrity measure}
\label{subsec: integrity measure}
The integrity measure provides a measure of the impact of external forcing on the stability of the ship.
For our purposes the integrity measure is defined as the Lebesgue measure of initial conditions within the region $U = \{(x,y,v_x,v_y): H(x,y,v_x,v_y) < \tfrac14 , -1 < y < 1 \}$ that do not lead to capsize. This region $U$ is contained within the basin of attraction of the origin for the autonomous system, and so any decrease corresponds to an erosion of the basin of attraction for the upright trajectory.  The region $U$ can be bounded by a rectangular region that we call $C$:
\begin{equation}
\begin{aligned}
    C=\{(x,y,v_x,v_y): -1<x<1, -1 < y < 1, -\sqrt{h} < v_x < \sqrt{h}, -\sqrt{2} < v_y < \sqrt{2}  \},
\end{aligned}
 \end{equation}
 which makes it easier to sample points. The volume of $U$ can be computed analytically by evaluating a definite integral: 
\begin{equation}
    \label{eq: volU}
    V_U = \int_{-1}^1 dy \int_{H(x,y,v_x,v_y) < \tfrac14}dxdv_xdv_y = \frac{163\pi \sqrt{2h}}{210} \approx 3.45 h.
\end{equation} 
Points are randomly selected within $C$ and if they lie within $U$ they are then classified as a capsize or non-capsize point by checking on which side of the stable manifold they lie. The probability of capsize is then evaluated by calculating the number of capsize points divided by the number of points within $U$.
Table \ref{tab: volume of safe region}  presents the integrity measure  of each of the systems given by the different forcing functions.
\begin{table}[hbtp]
    \centering
    \begin{tabular}{|c|c|c|}
    \hline
    \textbf{Forcing} & \textbf{Volume of safe} & \textbf{Relative volume}  \\
     \textbf{function} & \textbf{region in $U$} & \textbf{of safe region} \\
    \hline
    Autonomous & 3.45 & 1.00 \\
    Quasi-periodic forcing & 1.45 & 0.42\\
    Filtered white noise forcing & 1.38  &  0.40 \\
    \hline
    \end{tabular}
    \vspace*{5mm}
    \caption{Integrity measure for the system given by \eqref{eq: 2dof rescaled system}, with quasi-periodic and filtered white noise forcing. The relative volume is the volume of the safe region in $U$ divided by $V_U$.}
    \label{tab: volume of safe region}
\end{table}

\section{Discussion and conclusions}
\label{sec: discussion and conclusions}
The approach outlined in this paper  can be summarised as follows:
\begin{enumerate}
    \item We examined the autonomous dissipative system and linearised the non-linear vector field about the saddle points.
    \item For quasi-periodic and filtered white noise forcing we derived the hyperbolic trajectories for the linearised system. Through the algorithm in Section \ref{sec: algorithm hyp traj}, we found the hyperbolic trajectories for the non-linear system.
    \item Applying the algorithm in Section \ref{subsec: algorithm stable manifold} we computed a number of trajectories along the stable and centre manifold. These points were  used to reconstruct the stable manifold by interpolating the immersion $\iota_{S}$ given by this algorithm.
    \item By considering upon which side of the stable manifold an initial condition lies we determined the region of initial conditions which lead to transition over the saddle point/ centre manifold. We implemented this for a one degree of freedom system transition over an energy maximum and a two degree of freedom system modelling ship capsize. 
    \item For the latter system, we systematically computed the time to capsize using a dividing manifold, obtained using the stable manifolds. We also computed an integrity measure of the ship, to describe quantitatively how forcing affects the ship.
\end{enumerate}
The main significance of this work is that it allows aperiodic forcing to be taken into account in the analysis of escape from a potential well,  which is desired for many applications.
It is expected that the method can be extended without major difficulty to systems with more degrees of freedom.

Possible areas of extension are:
\begin{enumerate}
    \item Incorporate probability distributions over forcing functions and initial conditions to obtain probability distributions of the time to capsize over initial conditions.
    \item Extension to higher degrees of freedom. An energy based approach to modelling the 6 DoFs that includes configuration dependence of kinetic energy and potential energy and dissipative forcing due to waves was presented in \cite{nayfeh1974perturbation}. 
    \item The establishment of control strategies, to keep the ship on the correct side of the stable manifolds.
    \item Computation of a normally hyperbolic submanifold without requiring it to be a centre manifold for a hyperbolic trajectory.  We used dissipation in this paper to guarantee existence of a hyperbolic trajectory corresponding to the saddle of the unforced system, and then to find its centre manifold, but in principle we could have computed a normally hyperbolic submanifold directly without having to find a hyperbolic trajectory.  This will be important when extending to include neutral directions such as yaw, surge and sway.
\end{enumerate}

\backmatter

\bmhead{Supplementary material section} We do not provide supplementary material for this study.

\bmhead{Data availability statement}
The data that supports the findings of the study are openly available at \url{https://github.com/AlexMcSD/ship-capsize-2022}.
\bmhead{Acknowledgments}
This work was supported by the London Mathematical Society through a summer student bursary to the first author in 2022, with matched funding from the Mathematics Institute, University of Warwick.



\section*{Declarations}
{For the purpose of open access, the authors have applied a Creative Commons Attribution (CC-BY) licence to any Author Accepted Manuscript version arising from this submission.}
\bmhead{Conflict of interests statement}
The authors declare they have no competing interests.
\bmhead{Author contributions}
Alex McSweeney-Davis, MMath (Investigation: Lead; Software: Lead; Writing – original
draft: Lead; Writing - review and editing : Equal) \\
Robert Sinclair MacKay, PhD, FRS, FInstP, FIMA (Conceptualization: Lead; Funding
acquisition: Lead; Methodology: Lead; Supervision: Equal; Writing – original draft:
Supporting; Writing – review and editing: Equal) \\
Shibabrat Naik, PhD (Methodology: Supporting; Supervision: Equal; Writing – original
draft: Supporting; Writing – review and editing: Equal, Software: Supporting).

\noindent






\appendix
\section{Derivation of non-dimensionalised system}
\label{appendix: derivation}
The system given in Equation \eqref{eq: thompson  equations of motion} is a Hamiltonian system with Lagrangian given by:
\begin{align}
    \mathcal{L} &= \tfrac12 \left( I \dot{\phi}^2 + m \dot{z}^2\right) - V(z,\phi), \\
V(z,\phi)&= \frac{1}{2} c \phi_{v}^{2}\left[\left(\phi / \phi_{v}\right)^{2}-\frac{1}{2}\left(\phi / \phi_{v}\right)^{4}\right]+\frac{1}{2} h\left[z-\frac{1}{2} \gamma \phi^{2}\right]^{2}.
\end{align}
To this system, linear damping forces are added to each equation to yield
\begin{equation}
    \begin{pmatrix}
        \ddot{z} \\
        \ddot{\phi}
    \end{pmatrix} = -\nabla V(z,\phi) - \begin{pmatrix}
       k_z \dot{z} \\
       k_\phi \dot{\phi}
    \end{pmatrix}.
\end{equation}
To reduce the number of parameters, we introduce dimensionless quantities  $x = z/z_v,y = \phi/\phi_v, s = \sqrt{\frac{c}{I}}t$, where $z_v = \tfrac12 \gamma \phi_v^2$. Following \cite{thompson1996suppression} it is assumed that $\gamma^2 = \frac{2c}{h \phi_v^2}$, and under these new variables the equations of motion become 
\begin{equation}
    \label{eq: 2dof rescaled system}
    \begin{split}
        x'' +\tilde h (x - y^2) + k_x x' &= \tilde F (s), \\
y'' +y - xy + k_y y' &= \tilde M (s),
    \end{split}
\end{equation}
where $\tilde h =h\frac{I}{mc} $, $ k_x =  k_z\frac{\sqrt{c}}{m\sqrt{I}}$, $ k_y =  k_\phi \sqrt{\frac{1}{cI}}$, $ \tilde F (s) =  \frac{I}{m c z_v}F\left( \sqrt{\frac{I}{c}} s\right)$,   $ \tilde M (s) = \frac{1}{c \phi_v} M\left( \sqrt{\frac{I}{c}} s\right) $.
The Hamiltonian to the rescaled system is given by 
\begin{equation}
H(x,y,v_x,v_y) = \frac{1}{4 h}(v_x)^2 + \frac{1}{2}(v_y)^2 +\frac{1}{2}\left( y^2 + \frac{1}{2}x^2 - xy^2 \right),\end{equation}
and is a Lyapunov function in the unforced case.
In the main text of this paper only the rescaled system in Equation  \eqref{eq: 2dof rescaled system} will be studied, and all tildes will be dropped from the parameters and functions. It is then clear that this system is equivalent to the one from Equation \eqref{eq: 2dof rescaled vector}.

\section{Equilibria of damped autonomous system}
\label{appedix: analysis}
We analyse the autonomous form of the 2DoF rescaled system, from Equation \eqref{eq: 2dof rescaled vector}.
The linearisation of $\mathbf{f}$ about the saddle points $(1,\pm 1,0,0)$ is given by the matrix
\begin{equation}D\mathbf{f}(1, \pm 1,0,0 ) = \begin{pmatrix}
     0 & 0 & 1 &0 \\
     0 & 0 & 0 & 1 \\
     - h & \pm 2 h & -k_x & 0 \\
     \pm 1 & 0 & 0 & -k_{y,0}
\end{pmatrix}, \end{equation}
where $k_{y,0} =  D'(0) >0$. 
This has a characteristic polynomial at both saddle points of:
\begin{equation}
    \label{eq: char poly for saddles}
    \chi (\lambda) = \lambda^4 + (k_x + k_{y,0})\lambda^3 + (k_xk_{y,0} + h)\lambda^2 + hk_{y,0} \lambda- 2h .
\end{equation}
In general the roots to Equation  \eqref{eq: char poly for saddles} are best found numerically as the number of parameters are too high for the algebraic solutions to be practical. 

For the purposes of demonstrating the qualitative nature of the model we shall 
take $k_x = k_{y} = k >0$. The eigenvalues are then given by:
\begin{align} 
\lambda_1 &= \frac{1}{2}\left(-k-\sqrt{k^{2} - \alpha}\right), \\
\lambda_2 &= \frac{1}{2}\left(-k+\sqrt{k^{2} - \alpha}\right), \\
\lambda_3 &= \frac{1}{2}\left(-k-\sqrt{k^{2} + \beta}\right) ,\\
\lambda_4 &= \frac{1}{2}\left(-k+\sqrt{k^{2} + \beta}\right),
\end{align}
where $\alpha = 2 h+2 \sqrt{h} \sqrt{8+h}$ and $\beta = -2 h+2 \sqrt{h} \sqrt{8+h}$.
It follows that the equilibrium is a saddle with one expanding direction.
The eigenvectors for $D\mathbf{f}(1,1,0,0)$ are given by:
\begin{align}
\mathbf{v}_1 &= \begin{pmatrix}
    \frac{k}{2} - \frac{1}{2} \sqrt{k^2 - \alpha } \\
    - \frac{\sqrt{h}}{8 \sqrt{h}} \left(-k  + \sqrt{k^2 - \alpha}  \right) - \frac{\sqrt{8+h}}{8 \sqrt{h}} \left( k - \sqrt{k^2 - \alpha} \right)\\
    -\frac{h}{2} - \frac{1}{2} \sqrt{h}\sqrt{8+h}\\
    1
\end{pmatrix}, \\
\mathbf{v}_2 &= \begin{pmatrix}
    \frac{k}{2} + \frac{1}{2} \sqrt{k^2 - \alpha }\\
    - \frac{\sqrt{h}}{8 \sqrt{h}} \left(-k  -\sqrt{k^2 - \alpha}  \right) - \frac{\sqrt{8+h}}{8 \sqrt{h}} \left( k+\sqrt{k^2 - \alpha} \right) \\
   -\frac{h}{2} - \frac{1}{2} \sqrt{h}\sqrt{8+h} \\
    1
\end{pmatrix},  \\
\mathbf{v}_3 &= \begin{pmatrix}
    \frac{k}{2} - \frac{1}{2} \sqrt{k^2 +\beta }\\
   - \frac{\sqrt{h}}{8 \sqrt{h}} \left(-k + \sqrt{k^2 + \beta} +  \right)- \frac{\sqrt{8+h}}{8 \sqrt{h}} \left(-k +\sqrt{k^2 +\beta} \right) \\
   -\frac{h}{2} + \frac{1}{2} \sqrt{h}\sqrt{8+h} \\
    1
\end{pmatrix},  \\
\mathbf{v}_4 &= \begin{pmatrix}
    \frac{k}{2} + \frac{1}{2} \sqrt{k^2 +\beta }\\
    - \frac{\sqrt{h}}{8 \sqrt{h}} \left(-k -\sqrt{k^2 + \beta}  \right)- \frac{\sqrt{8+h}}{8 \sqrt{h}} \left(- k  - \sqrt{k^2 +\beta} \right)\\
    -\frac{h}{2} + \frac{1}{2} \sqrt{h}\sqrt{8+h}\\
    1
\end{pmatrix} .
\end{align}
The eigenvectors for $D\mathbf{f}(1,-1,0,0)$ can be found similarly.
This gives different phase portraits depending on the sign of $k^2 - \alpha$. In the case with no damping this would be negative, giving two eigenvalues on the imaginary axis. If the damping is small compared to the hydrostatic forces, we would expect that $k^2 - \alpha <0$, and we assume that this  holds for the remainder of the paper. In this case we can split $\mathbf{v}_1$ and $\mathbf{v}_2$ into the real and imaginary parts by:
\begin{equation} \mathbf{v}_1 = \mathbf{v}^{(1)} + i \mathbf{v}^{(2)}, \quad \mathbf{v}_2 = \mathbf{v}^{(1)} - i \mathbf{v}^{(2)}, \end{equation}
where
\begin{align}
    \mathbf{v}^{(1)} &= \left(\frac{k}{2}, - \frac{1}{8 \sqrt{h}} \left(-k \sqrt{h} + k \sqrt{8+h} \right),-\frac{h}{2} - \frac{1}{2} \sqrt{h}\sqrt{8+h},1 \right)^T, \\
\mathbf{v}^{(2)} &= \left( -\frac{1}{2} \sqrt{-k^2 + \alpha }, - \frac{1}{8 \sqrt{h}} \left( \sqrt{h}\sqrt{-k^2 + \alpha} - \sqrt{8+h}\sqrt{-k^2 + \alpha} \right),0,0 \right)^T.
\end{align}
In this case the homogeneous solution to the linearised IVP with initial values $\mathbf{x}_0 = X_1 \mathbf{v}^{(1)} + X_2 \mathbf{v}^{(2)} + X_3 \mathbf{v}_3 +X_4 \mathbf{v}_4$ is given by:
\begin{equation}
    \label{eq: 2dof general linear sol}
    \begin{split}
    \mathbf{x}(t) &= e^{-kt/2} \left( X_1 \cos (\omega t) + X_2 \sin (\omega t)) \right) \mathbf{v}^{(1)}\\
    &+ e^{-kt/2} \left( -X_1 \sin (\omega t) + X_2 \cos (\omega t)) \right) \mathbf{v}^{(2)} \\
    &+ e^{\lambda_3 t} X_3 \mathbf{v}_3 + e^{\lambda_4 t} X_4 \mathbf{v}_4 
    \end{split}
\end{equation}
where $\omega =\tfrac12 \sqrt{\alpha - k^2}$. We can transform this back into the original coordinates by substituting the values for the eigenvectors.

\section{Flux over the saddle in the autonomous undamped case}

\label{appendix:autoHam}

In the absence of forcing and damping, the dynamics simplify because energy is conserved.  So, the system can be studied on one energy surface at a time. Here, we apply the standard theory for transport across saddles of autonomous Hamiltonian systems to the (unforced, undamped) roll-heave model. 

At an energy above the saddle energy, $E_s$, that is $E > E_s$, there is a hyperbolic periodic orbit associated with the saddle. The forward and backward contracting invariant manifolds of the hyperbolic periodic orbit are diffeomorphic to $\mathbb{R}^1 \times \mathbb{S}$.
To compute these manifolds, first, we compute the hyperbolic periodic orbit at energy, $E > E_s$, using a numerical method called differential correction described in detail for this system in~\cite{naik2017geometry}. This method starts from a small amplitude ($\approx 10^{-4}$ for the system at hand) periodic solution (corresponding to the complex eigenvalues) to the linear system at the saddle and makes a small correction based on the state transition matrix to obtain the next periodic solution. This method gives a family of hyperbolic periodic orbits at energies above the saddle energy. Once the periodic orbit at the desired energy is crossed, a bisection approach is used to obtain the periodic orbit at the desired energy. Then, to compute the forward (backward) contracting manifold, we initialize points along the eigenvectors corresponding to the negative (positive) real eigenvalue and globalize them by integrating them backward (forward) with the full nonlinear equations of motion. 

\begin{figure}
    \centering
    \includegraphics[width = 0.8\textwidth]{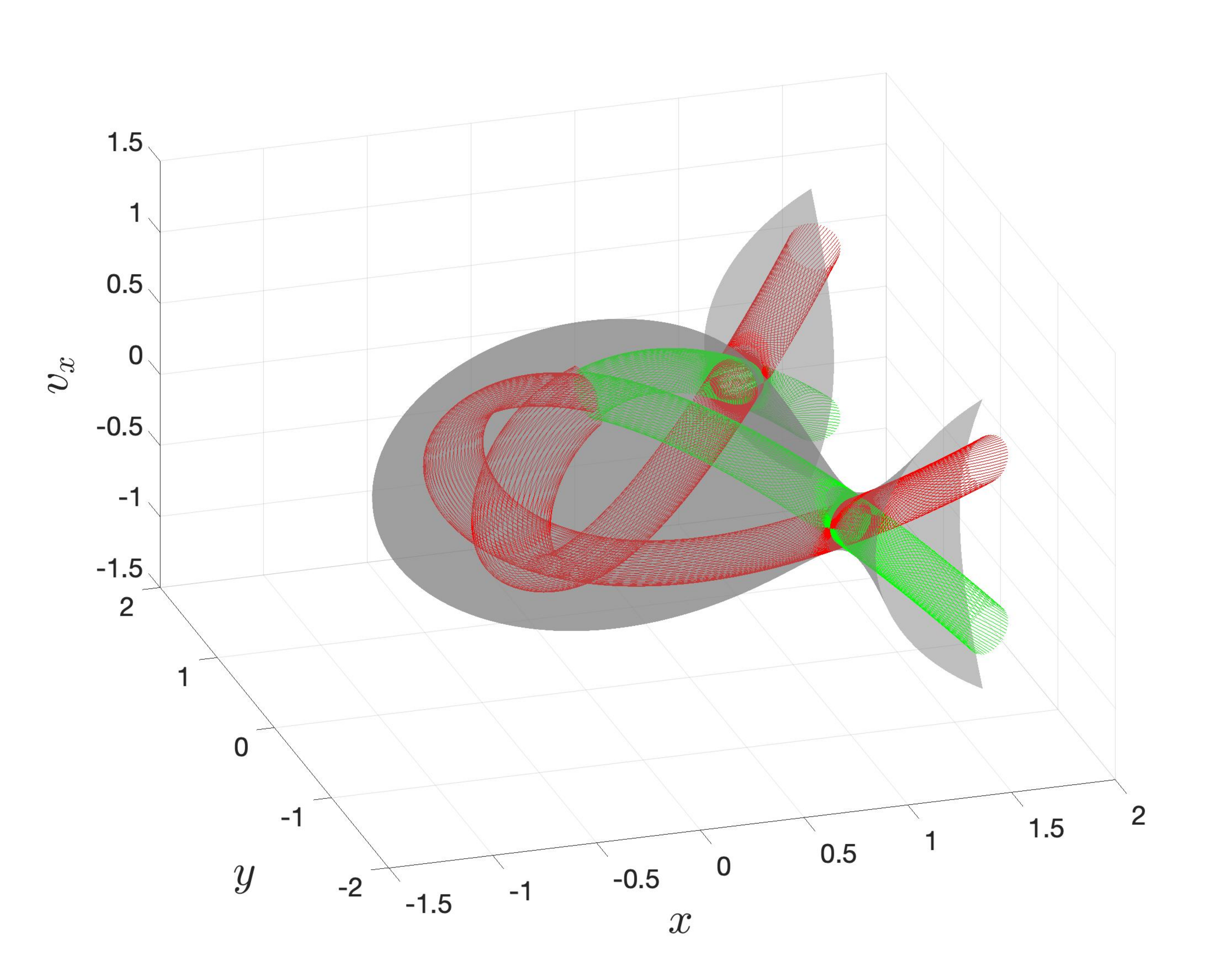}
    \caption{Forward (shown as green trajectories) and backward (shown as red trajectories) contracting invariant manifolds and the projection of the energy surface (shown as a grey surface) at $E = 0.26$. The manifolds are integrated until they intersect a section $\{ (x,y,v_x,v_y) \; \vert \; x = 0, v_x > 0 \}$ to aid in visualizing the manifolds' geometry.}
    \label{fig:auto_manifolds_3d}
\end{figure}

These manifolds are shown in Fig.~\ref{fig:auto_manifolds_3d} for $E = 0.26, h = 1$ (see~\cite{naik2017geometry} for more details on the numerical method). 
The figure shows their projection to $(x,y,v_x)$.  Also shown in grey is the surface where $v_y=0$ with the given energy. The projection of the energy surface is the region bounded by the grey surface, but the energy surface is a double covering of this region, given by two solutions for $v_y$ that merge at the boundary (although this is a simple way to represent an energy surface for a 2DoF system, it is better to use one-to-one representations, as in \cite{BMN} for example, but they require more work).
Since the contracting manifolds of the periodic orbits are codimension-1 in the 3D energy surface, they partition the volume into those points that will cross the saddle or remain inside the well around $(0,0)$. Furthermore, for each of the saddles, the intersections of the forward and backward manifolds with an appropriate 2D section will generate a homoclinic tangle and lobes, as has been studied by~\cite{mcrobie1991lobe}.
\bibliography{bibliography}

\end{document}